\documentclass[a4paper,11pt,leqno]{article}
\usepackage{amsmath,amsfonts,amsthm,amssymb}
\usepackage[left=3cm, right=3cm, top=3cm, bottom=3.5cm]{geometry}
\usepackage{graphicx}
\usepackage{mathrsfs}
\usepackage{tikz}
\usepackage{marginnote}
\numberwithin{equation}{section}
\usepackage[final]{pdfpages}

\theoremstyle{plain}
\newtheorem{theorem}{Theorem}[section]

\newtheorem{lemma}[theorem]{Lemma}
\newtheorem{proposition}[theorem]{Proposition}

\theoremstyle{remark}
\newtheorem{remark}{Remark}

\theoremstyle{definition}

\newcommand{\R}{\mathbb{R}}

\newcommand{\SF}{\mathbb{S}}
\newcommand{\eps}{\epsilon}

\usepackage{color}

\begin{document}

\title{Sharp lower bounds for the vector Allen-Cahn energy and qualitative properties of minimizers under no symmetry hypotheses}
\author{
{ Nicholas D. Alikakos\footnote{Department of Mathematics, University of Athens, Panepistemiopolis, 15784 Athens, Greece; e-mail:{\texttt{nalikako@math.uoa.gr}}}}$^{,}$\
{ Giorgio Fusco\footnote{Dipartimento di Matematica Pura ed Applicata, Universit\`a degli Studi dell'Aquila, Via Vetoio, 67010 Coppito, L'Aquila, Italy; e-mail:{\texttt{fusco@univaq.it}}}}
}
\date{}
\maketitle

\begin{abstract}
We study vector minimizers $u^\epsilon$ of $J^\epsilon_\Omega
(u) =\int_R(\frac{\epsilon}{2}\vert\nabla u\vert^2+\frac{1}{\epsilon}W(u))dx$, $W>0$ on $\R^m\setminus\{a_1,\ldots,a_N\}$, $m\geq 1$ for bounded domains $\Omega\subset\R^n$
, with certain geometrical features
and $u = g_\epsilon$ on $\partial\Omega$
. We derive a sharp lower bound (as $\epsilon\rightarrow 0$) with the additional feature
that it involves half of the gradient and part of the domain
. Based on this we derive very
precise (in $\epsilon$) pointwise estimates up to the boundary for $\lim_{\epsilon\rightarrow0} u^\epsilon = u^0$. Depending on
the geometry of $\Omega$
 $u^\epsilon$ exhibits either boundary layers of internal layers. We do not impose
symmetry hypotheses and we do not employ $\Gamma$-convergence techniques.
\end{abstract}
\section{Introduction}
The object of study in the present paper is the system
\begin{equation}
\Delta u-W_u(u)=0,\;\;x\in\Omega\subset\R^n,
\label{elliptic}
\end{equation}
where
 $W:\R^m\rightarrow\R$ is a $C^2$ phase transition potential. That is: $W$ is  nonnegative and vanishes only on a finite set $\{W=0\}=A=\{a_1,\ldots,a_N\}$ for some distinct points $a_1,\ldots,a_N\in\R^m$ that represent the phases of a substance that can exist in $N\geq 2$ different equally preferred states. We assume that the zeros $a_1,\ldots,a_N$ are nondegenerate in the sense that the jacobian matrix $W_{zz}(a)$, $a\in A$ is definite positive. Finally we assume that
\begin{equation}
\label{infty-bound}
\liminf_{\vert z\vert\rightarrow+\infty}W(z)>0.
\end{equation}
System \eqref{elliptic} is the Euler-Lagrange equation corresponding the the Allen-Cahn energy
\begin{equation}
J_\Omega(v)=\int_\Omega\Big(\frac{1}{2}\vert\nabla v\vert^2+W(v)\Big)dx.
\label{ACE}
\end{equation}
We are interested in the class of solutions that connect in some way the phases or a subset of them.

The scalar case $m=1$ has been extensively studied: here $N=2$ is the natural choice. The reader may consult \cite{Wei},\cite{sav}, \cite{dkw} where further references can be found. A well known conjecture of De Giorgi (1978) and its solution about thirty years later, played a significant role in the development of a large part of this work.

The vector case $m\geq 2$ by comparison has been studied very little. We note that for coexistence of three or more phases a vector order parameter is necessary and so there is physical interest for the system. On the geometric side
\eqref{ACE}, when rescaled as in \eqref{ACE-res} below, produces minimal partitions $\{P_j\}_1^N$ (see \eqref{perimeter}) that exhibit {\em Junctions} that is singularities with certain structures, that do not exist for $m=1$.

For $m\geq 2$ \eqref{elliptic} has been mainly studied in the class of equivariant solutions with respect to reflection groups beginning with Bronsard, Gui and Schatzman \cite{bgs} and later Gui and Schatzman \cite{gs}, and significantly extended and generalized in various ways by the authors including also \cite{BfS}. We refer also to Chapters 6 and 7 in \cite{afs} and the reference there in. The only related work that does not require symmetry of which we are aware of is Schatzman \cite{scha}.

We will be focusing on the rescaled functional
\begin{equation}
J_\Omega^\eps(v)=\int_\Omega\Big(\frac{\epsilon}{2}\vert\nabla v\vert^2+\frac{1}{\epsilon}W(v)\Big)dx,
\label{ACE-res}
\end{equation}
where $\Omega$ is open, bounded, $C^{1,\alpha}$, for some $\alpha\in(0,1)$, smooth  connected set and we consider the minimization problem
\begin{equation}
\min J_\Omega^\epsilon(v),\;\;v=g_\epsilon\;\;\text{on}\;\partial\Omega,
\label{Min_Prob}
\end{equation}
where $g_\epsilon$ is a given map that may depend on $\epsilon$.

The rescaled problem \eqref{Min_Prob} is also useful for constructing entire solutions of \eqref{elliptic} over $\R^n$, and although this point of view is not exploited in the present paper, it provides one of the motivations behind this work.

We are interested in uniformly (with respect to $\epsilon$) pointwise bounded global minimizers connecting minima of $W$, and for this reason we adopt the simple hypothesis

\begin{equation}
\begin{split}
& W_u(u)\cdot u> 0,\;\;\text{for}\;\vert u\vert> M,\;\text{some}\;M,\\
& \vert g_\epsilon\vert\leq M.
\end{split}
\label{u>M}
\end{equation}

A major general tool for the study of the minimizers $u_\epsilon$ of \eqref{Min_Prob} is the limiting functional as $\epsilon\rightarrow 0$ given by the weighted perimeter functional
\begin{equation}
E(\mathcal{P})=\sum_{i\neq j}^N\sigma_{i,j}\mathcal{H}^{n-1}(\partial D_i\cap\partial D_j),\;\;
\mathcal{P}=\{D_j\}_{j=1}^N,\;\;\text{partition of}\;\Omega.
\label{perimeter}
\end{equation}
Baldo \cite{B} established for the vector case $m\geq 2$ for the mass constrained problem, that $\Gamma-\lim_{\epsilon\rightarrow 0}J_\Omega^\epsilon=E(\mathcal{P})$, from which it follows that along subsequences \begin{equation}
\lim_{\epsilon\rightarrow 0}\Vert u_\epsilon-u_0\Vert_{L^1(\Omega;\R^m)}=0,
\label{L1-lim}
\end{equation}
where $u_0=\sum_{j=1}^Na_j{\chi}_{D_j}$ is a minimizer of $E$ under the same constraint.

The other major tool for the study of the minimizers $\{u_\epsilon\}$ as $\epsilon\rightarrow 0$ is the Caffarelli-Cordoba \cite{CC} density estimate originally derived for $m=1$. This is independent from $\Gamma$-convergence and complements it by upgrading \eqref{L1-lim} to uniform convergence over compacts in $D_j$.
 It has been extended to the vector case in \cite{A-F}. We refer to \cite{afs} Chapter 5.

A major difficulty that one faces in implementing these general results and their variants is that the convergence in \eqref{L1-lim} does not come with an estimate in $\epsilon$. This also is the major obstruction for utilizing the rescaled problem in constructing entire solutions to \eqref{elliptic}. It is only under symmetry conditions that such estimates have been obtained in generality. Our approach here produces such estimates for the problem at hand.

In the present paper we illustrate in term of various simple examples how the derivation of sharp lower bounds for $J_\Omega^\epsilon(u_\epsilon)$ allows extraction of pointwise estimates for $u_\epsilon\rightarrow u_0$ all the way up to the boundary of the partition. We do not make any symmetry assumptions. We do not utilize the limiting problem, but we do utilize the vector Caffarelli-Cordoba density estimate.

We now describe the content of the paper in more detail. We consider two examples where beside an upper bound it is also possible to derive a sharp lower bound for the energy of a minimizer of problem \eqref{Min_Prob}. We show that the knowledge of sharp lower and upper bounds together with the vector Caffarelli-Cordoba density estimate \cite{CC}, \cite{A-F} allow for an accurate description of the fine structure of minimizers for all $\epsilon>0$ sufficiently small. A key point of the analysis is that the lower bound is obtained by considering only the energy in a proper subset of $\Omega$.

In our first example $\Omega$ is a bounded domain in $\R^n$, $n\geq 2$ and $g_\epsilon\equiv z$ where $z\in\R^m\setminus A$ is a fixed vector. We prove
\begin{theorem}
Let $u$ ($u=u_\epsilon$) a minimizer of problem \eqref{Min_Prob} with $g_\epsilon\equiv z$, $z\in\R^m\setminus A$. Then there exist $a_+\in A$ and positive constants $k,K$ and $C$ such that
\[\vert u(x)-a_+\vert\leq Ke^{-\frac{k}{\epsilon}(d(x,\partial\Omega)-C\epsilon^\frac{1}{3(n-1)})^+},\;\;x\in\Omega.\]
\end{theorem}
The proof of this theorem is based on the fact that one can show (see Theorem \ref{lower}) that most of the energy of a minimizer is contained in a tiny neighborhood of $\partial\Omega$. This allows for a transparent use of the density estimate in combination with the sharp upper and lower energy bounds.

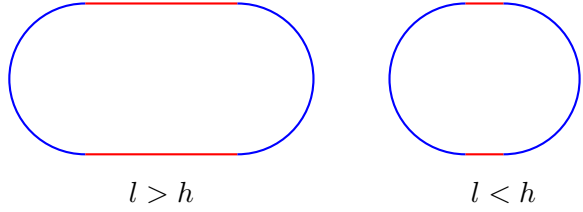
\begin{figure}
  \begin{center}
\begin{tikzpicture}
\draw [blue,thick] (0,1) arc [radius=1, start angle=90, end angle=270];
\draw [red,thick](0,1)--(2,1);
\draw [red,thick](0,-1)--(2,-1);
\draw [blue,thick] (2,-1) arc [radius=1, start angle=-90, end angle= 90];
\node at (1,-1.5){$l>h$};

\draw [blue,thick] (5,1) arc [radius=1, start angle=90, end angle=270];
\draw [red,thick](5,1)--(5.5,1);
\draw [red,thick](5,-1)--(5.5,-1);
\draw [blue,thick] (5.5,-1) arc [radius=1, start angle=-90, end angle= 90];
\node at (5.5,-1.5){$l<h$};
\end{tikzpicture}
\end{center}
\caption{blue$=a_-$, red$=a_+$. $l>h$ or $l<h$.}
\label{fig3}
\end{figure}

In our second example we consider a domain $\Omega\subset\R^2$, otherwise as introduced after \eqref{ACE-res} and with the additional geometrical property

\[
\Omega\cap[0,l]\times[0,h]=[0,l]\times(0,h).
\]
The model domain is as in Figure \ref{fig3}, although our definition allows for more general domains (see Figure \ref{OOmega}). Here $l>0$ and $h>0$ are free parameters and our interest is to describe the fine structure of minimizers as function of the ratio $h/l$.
We choose the Dirichlet data $u=g_\epsilon$ on $\partial\Omega$, $g_\epsilon$ a $C^{1,\alpha}(\bar{\Omega};\R^m)$, some $\alpha\in(0,1)$, with the feature that $g_\epsilon$ converges, as $\epsilon\rightarrow 0$, in a controlled manner to a step map taking  values $a_-, a_+$ in $A$. Here $a_-\in A$ can be fixed arbitrarily while $a_+$ is chosen via a minimization process (see Lemma \ref{E-conn}) that determines a minimal connection $\bar{u}:\R\rightarrow\R^m$ between $a_-$ and $a_+$. We assume

\[
\begin{split}
&g_\epsilon(x,0)=g_\epsilon(x,h)=a_+,\;\;x\in(C_0\epsilon,l-C_0\epsilon),\\
&g_\epsilon(x,y)=a_-,\;\;(x,y)\in\partial\Omega\setminus(0,l)\times\{0,h\},
\\
&\vert g_\epsilon(x,y)\vert\leq C,\;\vert g_{\epsilon,x}(x,y)\vert\leq \frac{C}{\epsilon}
,\;\;x\in(0,C_0\epsilon)\cup(l-C_0\epsilon,l),\;\;y=0,h\\
&\vert g_{\epsilon,x}(\cdot,y)\vert_\alpha\leq \frac{C}{\epsilon^{1+\alpha}},\;\;y=0,h.
\end{split}
\]

We show that the structure of global minimizers $u_\epsilon$ for $0<\epsilon<<1$ depends drastically on whether $l<h$ of $l>h$ that is on the geometry of $\Omega$.

We are stating somewhat loosely our main results for convenience of the reader

\begin{theorem}
\label{1-T}
$(l<h$, The Boundary Layer Case$)$
There is $\epsilon_0>0$ such that, if  $u_\epsilon$, $\epsilon\in(0,\epsilon_0]$ is a minimizer of \eqref{Min_Prob}, then
\begin{equation}
\begin{split}
&2\sigma l-C\epsilon\leq\int_R\Big(\frac{\epsilon}{2}\vert\frac{\partial u_\epsilon}{\partial y}\vert^2+\frac{1}{\epsilon}W(u_\epsilon)\Big)dx dy\leq J_\Omega^\epsilon(u_\epsilon)\leq2\sigma l+C\epsilon\vert\ln{\epsilon}\vert^3,\\\\
&\vert u_\epsilon(z)-a_-\vert\leq K e^{-\frac{k}{\epsilon}(d(z,\partial^+\Omega)-C\epsilon^\frac{1}{4}\vert\ln{\epsilon}\vert^\frac{3}{2})^+},\;\;z\in\Omega,
\end{split}
\label{l>hEst}
\end{equation}
where $R=(0,l)\times(0,h)$, $\sigma$ is the Action of the connecting orbit between $a_-$ and $a_+$ and $C$, $k$ and $K$ are positive constants, $\xi^+=\max\{0,\xi\}$ and $\partial^+\Omega=(0,l)\times\{0,h\}$.
\end{theorem}

These estimates imply that $u_\epsilon$ converges uniformly in compacts in $\bar{\Omega}\setminus\partial^+\Omega$ to $a_-$ and there is a boundary layer in a neighborhood of $\partial^+\Omega$ which can be shown to be strictly thicker that $\mathrm{O}(\epsilon)$ ($\mathrm{O}(\epsilon)/\mathrm{o}(1)$), see Theorem \ref{thick} below.

\begin{theorem}
\label{2-T}$(l>h$, The Internal Layer Case$)$
There is $\epsilon_0>0$ such that, if  $u_\epsilon$, $\epsilon\in(0,\epsilon_0]$ is a minimizer of \eqref{Min_Prob}, then
\begin{equation}\begin{split}
&2\sigma l-C\epsilon\leq\int_D\Big(\frac{\epsilon}{2}\vert\frac{\partial u_\epsilon}{\partial x}\vert^2
+\frac{1}{\epsilon}W(u_\epsilon)\Big)dx dy\leq J_\Omega^\epsilon(u_\epsilon)\leq2\sigma l+C\epsilon,\\\\
&\vert u_\epsilon(z)-a_-\vert\leq K e^{-\frac{k}{\epsilon}(d(z,R)-C\epsilon^\frac{1}{4})^+},\;\;z\in\Omega\setminus R,\\
&\vert u_\epsilon(z)-a_+\vert\leq K e^{-\frac{k}{\epsilon}(d(z,\Omega\setminus R)-C\epsilon^\frac{1}{4})^+},\;\;z\in R,
\end{split}
\label{l<hEst}
\end{equation}
where $D$ is a strict subset of $\Omega$ (cfr. Figure \ref{D}), $R=(0,l)\times(0,h)$,  $K$, $k$ and $C$ positive constants.
\end{theorem}
These estimates imply that, as $\epsilon\rightarrow 0$, $u_\epsilon$ converges to the step map $u_0$:
\[u_0=\left\{\begin{array}{l}
a_-,\;\;\text{in}\;\Omega\setminus R,\\
a^+,\;\;\text{in}\; R,
\end{array}\right.
\]
and the convergence is uniform in compacts in $\Omega\setminus\{0,l\}\times[0,h]$ and all the way up to $\partial\Omega$. We expect the internal layer to be of $\mathrm{O}(\epsilon)$ thickness with an one dimensional profile given by a heteroclinic connection.
In Figure \ref{fig3} below we exhibit for a simple geometry the two different situations.

At a formal level is not difficult to understand why minimizers favour the boundary layer or the internal layer by considering the limiting problem. Assume that the Dirichlet data are preserved and compare the interface energies. 
 We refer to \cite{ABP} where a scalar problem with $\epsilon$-dependent Dirichlet data similar to ours is considered and the $\Gamma$-limit is justified rigorously.

The lower bound estimates have different features that play different roles in the derivation of the pointwise estimates.

Fist note that they involve only part of the gradient. For instance in the internal layer case via the lower and upper bounds we obtain
\begin{equation}
\int_\Omega\vert\frac{\partial u_\epsilon}{\partial y}\vert^2dx dy\leq C,
\label{Kin-y}
\end{equation}
which implies that, for small $\epsilon>0$, the interfaces are almost orthogonal to the $x$ axis. We note that this point is reminiscent of an estimate in Alama Bronsard Gui \cite{ABG} and also in Schatzman \cite{scha} where however the setup and the arguments are entirely different (cfr. \cite{afs} Theorem 8.5 and Lemma 9.4). In the boundary layer case the presence of the logaritm term (that we don't think can be removed) does not allow to derive a bound independent of $\epsilon$ for $\int_\Omega\vert\frac{\partial u_\epsilon}{\partial x}\vert^2dx dy$. However it serves a similar purpose (Lemma \ref{Kin-bound}) for establishing the width of the boundary layer as we explain later in this introduction.

The other feature of the lower bounds is that they are based on a strict subset of $\Omega$. This fact, in conjunction with the upper bound can be used to obtain certain refinements and estimates up to the boundary in the remaining part of $\Omega$ via the density estimate.

The existence of the boundary layer in Theorem \ref{1-T} is a higher dimensional effect, and it is a new phenomenon. It is the outcome of a compromise between two competing minimization requirements. On the one hand, in the interior of $\partial^+\Omega$, $u_\epsilon$ tries to behave like in the one-dimensional case and pushes the layer in the interior. On the other hand in a neighborhood of the extreme points of $\partial^+\Omega$ minimization requires that the solution must remain near $a_-$. At a formal level it is not difficult to guess that the layer is strictly thicker than $\mathrm{O}(\epsilon)$ since there is no connecting orbit in the half line. For establishing this one needs to show that the problem in the half space

\begin{equation}
\begin{split}
&\Delta U^0=W_u(U^0),\;\;(\xi,\eta)\in\R\times(0,\infty),\\
&U^0(\xi,0)=a_+,\;\;\xi\in\R,
\end{split}
\label{eq-half}
\end{equation}
has the unique solutions $U^0\equiv a_+$. In the scalar case $m=1$, and in the whole plane (on $\R^n$ generally) it was established by Modica \cite{Mo} that, if the entire solution $U$ has a point $x_0\in\R^n$ such that $W(U(x_0))=0$, then necessarily $U\equiv a$ where $a=U(x_0)\in A$. He obtained this as a corollary of the so called Modica Inequality
\[\frac{1}{2}\vert\nabla U\vert^2\leq W(U),\]
that is valid for $L^\infty$ solutions of $\Delta U=W_u(U)$, $U:\R^n\rightarrow\R$.

Farina and Valdinoci \cite{FaVa}  have extended  the Modica inequality to solutions satisfying $\Delta U=W_u(U)$ in the half space $\R^{n-1}\times[0,\infty)$ provided that $U\in C^2(\R^{n-1}\times[0,\infty))\cap L^\infty$. Utilizing this result we settle \eqref{eq-half} for $m=1$.

However for $m\geq 2$ P.Smyrnelis \cite{S-P}  (cfr. 3.3 in \cite {afs}) has pointed out that the Modica inequality is not generally valid and thus a different approach is required. For this purpose we can employ an appropriate Hamiltonian Identity of Gui \cite{Gui}  (e.g. 3.4 in \cite{afs}) together with the analog of \eqref{Kin-y} mentioned above to deduce that
\begin{equation}
\int_\R\vert U_n^0(\xi,0)\vert^2d \xi=0,
\label{U0-nabal}
\end{equation}
and thus  conclude via a unique continuation Theorem in \cite{gl}  that $U^0(\xi,\eta)=a_+$.

The following simple estimate on the $1$-dimensional energy is basic and is implemented in several places and forms throughout the paper
\begin{equation}
J_{(s_1,s_2)}(v)=\int_{s_1}^{s_2}\Big(\frac{\vert v^\prime\vert^2}{2}+W(v)\Big)ds\geq \sigma-\frac{1}{2}(\delta_-^2+\delta_+^2),
\label{1-E}
\end{equation}
where $C_W$ is a constant and $\delta_\pm$ satisfy $\vert v(s_1)-a_-\vert\leq\delta_-$ and $\vert v(s_2)-a_+\vert\leq\delta_+$ (cfr. Lemma \ref{lower-sigma} below). The estimate \eqref{1-E} is, in particular utilized in the derivation of the lower bound.

The estimates \eqref{l<hEst} and \eqref{l>hEst} are obtained via linear elliptic theory from an estimate of the type
\begin{equation}
\vert u_\epsilon(z)-a\vert\leq\delta,\;\;z\in V\subset\Omega,\;\text{some open}\;V,
\label{type}
\end{equation}
where $\delta>0$ is a small number.

A route for obtaining \eqref{type} in a neighborhood of the part of $\partial\Omega$ where $u=a$ is by constructing a $V$ with $\partial V$ partially coinciding with $\partial\Omega$ and such that
\[\vert u_\epsilon(z)-a\vert\leq\delta,\;\;z\in \partial V.\]
Then we conclude via Theorem 4.1 in \cite{afs}. It should be noted that this argument utilizes the fact that $\Omega$ is $2d$ and $\partial V$ is a curve along which \eqref{1-E} can be implemented.

Finally a few comments on possible extensions of our results are in order. We have studied the case $l<h$ and $l>h$ and established a significant qualitative difference between the global minimizers. We expect that both type of solutions exist independently of the relation between $l$ and $h$ but as local minimizers. The case $l=h$ is open.
Here we expect two type of global minimizers. Extending the results to genuinely higher dimensional examples is not straightforward. A minimizer of the type of Theorem \ref{1-T} should require that the boundary $\partial\Omega$ is a minimal surface. The lower bound also should be significantly harder since the sections will not be one-dimensional.

The remaining part of the paper is structured as follows: in $\S 2$ we present the basic lemmas. In $\S 3$ we present a different Example 1.
In $\S 4$ we study the example covered in Theorems \ref{1-T} and \ref{2-T} above. First we introduce the hypothesis in $\S 4.1$. Then in $\S 4.2$ we consider the boundary layer case and in $\S 4.3$ the internal layer case.

\section*{Acknowledgement}
We would like to thank V.Stefanopoulos for his comments on a previous manuscript that led to an improved presentation of $\S 2$.

\section{Basic Lemmas}\label{1}

\begin{lemma}
\label{Wlemma}
The assumptions on $W$ imply the existence of $\delta_W>0$, and constants $c_W, C_W>0$ such that
\begin{equation}
\begin{split}
&\vert u-a\vert=\delta,\\
&\Rightarrow\;\;\frac{1}{2}c_W^2\delta^2\leq W(u)\leq\frac{1}{2}C_W^2\delta^2;\;\delta\leq\delta_W,\;a\in A.
\end{split}
\label{W}
\end{equation}
Moreover, if  $\delta_W>0$ is sufficiently small, it results
\begin{equation}
\begin{split}
&\delta\in(0,\delta_W] \:\text{and}\;\min_{a\in A}\vert u-a\vert\geq\delta,\\
&\Rightarrow\;\;\frac{1}{2}c_W^2\delta^2\leq W(u).
\end{split}
\label{outW}
\end{equation}
\end{lemma}
For a map $v:(s_1,s_2)\rightarrow\R^m$ in $H_{\mathrm{loc}}^1$, with $-\infty\leq s_1<s_2\leq+\infty$, we define
\[J_{(s_1,s_2)}(v)=\int_{s_1}^{s_2}\Big(\frac{\vert v^\prime\vert^2}{2}+W(v)\Big)ds.\]
For $J_{(-\infty,+\infty)}(v)$ we also use the notation $J_\R(v)$.
\begin{lemma}
\label{E-conn}
Given $a_-\in A$ there exists a map $\bar{u}:\R\rightarrow\R^m$ that satisfies
\[\begin{split}
&\lim_{s\rightarrow-\infty}\bar{u}(s)=a_-,\quad\quad\lim_{s\rightarrow+\infty}\bar{u}(s)=a_+\in A\setminus\{a\},\\
&J_\R(\bar{u})=\min J_\R(\bar{v}),
\end{split}\]
where the minimization is taken on the set of $v\in H_{\mathrm{loc}}^1(\R;\R^m)$ that satisfy
\[\lim_{s\rightarrow-\infty}v(s)=a,\quad\quad \lim_{s\rightarrow+\infty}d(v(s), A^\prime)=0.\]
\vskip.2cm
Given $z\in\R^m\setminus A$ there exists a map $\bar{u}_+\in H_{\mathrm{loc}}^1(0,+\infty);\R^m)$ that satisfies
\[\begin{split}
&\lim_{s\rightarrow 0_+}\bar{u}_+(s)=z,\quad\quad\lim_{s\rightarrow+\infty}\bar{u}_+(s)=a_+\in A,\\
&J_{(0,+\infty)}(\bar{u}_+)=\min J_{(0,+\infty)}(v),
\end{split}\]
where the minimization is taken on the set of $v\in H_{\mathrm{loc}}^1(0,+\infty);\R^m)$ that satisfy
\[v(0)=z,\quad\quad \lim_{s\rightarrow+\infty}d(v(s), A)=0.\]

\end{lemma}
\begin{proof}
  See for example Theorem 2.1 in \cite{afs}) and its proof.
\end{proof}
\vskip.2cm
We set $\sigma=J_\R(\bar{u})$ and $\sigma^+=J_{(0,+\infty)}(\bar{u}_+)$. That is: $\sigma$ is the energy of the map $\bar{u}$ that connects $a$ to $a_+$ and $\sigma^+$ is the energy of the map $\bar{u}_+$ that connects $z$ to $a_+$.

Set  $\Gamma_0(a_\pm)=\{a_\pm\}$ and $\Gamma_\delta(a_\pm)=\partial B_\delta(a_\pm)$ for $\delta>0$.

\begin{lemma}
\label{lower-sigma}
Let $a_\pm$ ,$\bar{u}$ and $\bar{u}_+$ as in Lemma \ref{E-conn}.
\begin{enumerate}
\item Let $\delta_-\geq 0$ and $\delta_+\geq 0$ small numbers and
let $v:(s_-,s_+)\rightarrow\R^m$ a smooth map such that
\begin{equation}
\lim_{s\rightarrow s_\pm}d(v(s),\Gamma_{\delta_\pm}(a_\pm))=0.
\label{Lim-}
\end{equation}

Then
\begin{equation}
J_{(s_-,s_+)}(v)\geq\sigma-\frac{1}{2}C_W(\delta_-^2+\delta_+^2).
\label{sigmadelta}
\end{equation}
\item
let $v:(0,s_+)\rightarrow\R^m$ a smooth map such that
\begin{equation}
\begin{split}
&\lim_{s\rightarrow 0_+}v(s)=z,\\
&\lim_{s\rightarrow s_+}d(v(s),\Gamma_{\delta_+}(a_+))=0.
\end{split}
\label{Lim-1}
\end{equation}
Then
\begin{equation}
J_{(0,s_+)}(v)\geq\sigma^+-\frac{1}{2}C_W\delta_+^2.
\label{sigmadelta+}
\end{equation}
\end{enumerate}
\end{lemma}
\begin{proof}
1. For $\delta_-=\delta_+=0$ \eqref{sigmadelta} is just the statement of the minimality of $\bar{u}$. Therefore we can assume that either $\delta_-$ or $\delta_+$ or both are positive. From \eqref{Lim-} and $\delta_+>0$, if $s_+=+\infty$, it follows $\int_{s_-}^{s_+}W(v)ds=+\infty$ and \eqref{sigmadelta} holds trivially. The same is true if $\delta_+>0$, $s_+<+\infty$ and $\lim_{s\rightarrow s_+}v(s)$ does not exist. Indeed in this case we have $\int_{s_-}^{s_+}\vert\dot{v}\vert^2 ds=+\infty$. It follows that, if $\delta_+>0$, we can assume $s_+<+\infty$ and moreover that
\begin{equation}
\lim_{s\rightarrow s_+}v(s)=v_+,
\label{Lim}
\end{equation}
for some $v_+\in\Gamma_{\delta_+}(a_+)$. Analogous conclusion applies to the case $\delta_->0$.

2. If both $\delta_-$ and $\delta_+$ are positive and $w_\pm$ is a test map that connects $v_\pm$ to $a_\pm$, the minimality of $\bar{u}$ implies
\[J_{(s_-,s_+)}(v)\geq\sigma-J(w_-)-J(w_+),\]
where $J(w_\pm)$ is the energy of $w_\pm$. This yields \eqref{sigmadelta} provided we show that $w_\pm$ can be chosen so that
\[J(w_\pm)\leq\frac{1}{2} C_W\delta_\pm^2.\]

3. We choose
\[w_+=(1-\frac{\gamma(s)}{\delta_+})a_+ +\frac{\gamma(s)}{\delta_+}v_+,\;\;\gamma(s)=\delta_+e^{-C_W(s-s_+)}.\]
it follows, using also \eqref{W}

\[\begin{split}
&\frac{1}{2}\int_{s_+}^{+\infty}\vert\dot{w}_+\vert^2 ds=\frac{C_W^2}{2}\vert v_+-a_+\vert^2
\int_{s_+}^{+\infty}e^{-2C_W(s-s_+)}ds=\frac{1}{4}C_W\delta_+^2,\\
&\int_{s_+}^{+\infty}W(w_+)ds\leq\frac{C_W^2}{2}\vert v_+-a_+\vert^2\int_{s_+}^{+\infty}e^{-2C_W(s-s_+)}ds=\frac{1}{4}C_W\delta_+^2,
\end{split}\]
 This and the analogous computation for $J(w_-)$ establish \eqref{sigmadelta} for $\delta_-$ and $\delta_+$ positive. Clearly \eqref{sigmadelta} is valid also if $\delta_-$ or $\delta_+$ vanishes. The proof of
 \eqref{sigmadelta+} is analogous. The proof is complete.
\end{proof}

 We also define $\sigma^*$ by setting
\[
\sigma^*=\inf_vJ_{(s_1,s_2)}(v),
\]
where  $v\in H_{\mathrm{loc}}^1((s_1,s_2);\R^m)$ is a map that satisfies $v(s_1)=z$ and $\lim_{s\rightarrow s_2}d(v(s_2),A\setminus\{a_+\})=0$.

 Note that, if the extreme $a_+\in A$ of $\bar{u}_+$ is uniquely determined, then we have
\[\sigma^+<\sigma^*.\]
The same is true if $a_+$ is not unique provided the set of $v$ is restricted to the maps that satisfy
$\lim_{s\rightarrow s_2}d(v(s_2),A\setminus\{a\in A:a= a_+\})=0$.

\begin{figure}
  \begin{center}
\begin{tikzpicture}
\draw [](-5,-1.5).. controls (-4.5,-1)..(-4,.7);
\node at (-4.5,0){$S$};

\draw [](-1.8,1).. controls (-1,-.5) ..(1,-2);
\node at (-1.4,1.2){$S^\prime$};

\draw [thick,->](-4.45,-.71)--(-3.8,-.69);
\draw [thick,->](-4.45,-.71)--(-3.9,-.9);
\draw [thick,->](-.7,-.65)--(-.3,-.25);

\node at (-3.9,-.43){$e$};

\node at (-4,-1.25){$\nu$};

\node at (-.3,-.05){$\nu^\prime$};

\draw (-7,-1)--(1,-1);
\draw (-7,-.8)--(1,-.6);
\draw (-7,-.6)--(1,-.1);
\node at (-4.7,-.72){$d S$};
\node at (-3.2,-.55){$\omega$};


\end{tikzpicture}
\end{center}
\caption{The region $\omega$.}
\label{fig1}
\end{figure}
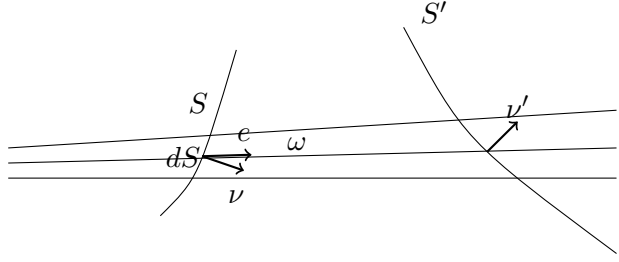

\begin{lemma}
\label{basic}
Consider a smooth family of lines that are transversal to two distinct $n-1$ dimensional surfaces $S$ and  $S^\prime$. Consider a point $p\in S$ and let $e$ be a unit vector parallel to the line of the family through $p$. Let $\nu$ a unit vector normal to $S$ at $p$. Let $dS$ a small neighborhood of $p$ in $S$ and let $dS^\prime$ the set of the intersections of the lines through $dS$ with $S^\prime$. Finally let $\omega$ the union of all segments determined on the lines of the family by $dS$ and $dS^\prime$.
Let $v:O\subset\R^n\rightarrow\R^m$, $O\supset\omega$ open, a smooth map that satisfies
\[\begin{split}
&v=z,\;\;x\in S,\\
&\vert v-a\vert\leq\delta,\;\;x\in S^\prime.
\end{split}\]
Then
\[J_\omega(v)\geq \sigma_\delta\min\{\nu\cdot e dS,\nu^\prime\cdot e dS^\prime\}.\]
This inequality still holds if the condition $\vert v-a\vert\leq\delta,\;\;x\in S^\prime$ is replaced by the condition that each segment in $\omega$ contains a point where $\vert v-a\vert\leq\delta$.
\end{lemma}
\begin{proof}
Follows from Fig. \ref{fig1}.
\end{proof}

\section{Example 1}\label{ex1}

 We consider the problem
 \begin{equation}
 \begin{split}
 &\min_{u\in\mathscr{A}} J_\Omega^\epsilon(u),\;\;J_\Omega^\epsilon(u)=\int_\Omega\Big(\frac{\epsilon}{2}\vert\nabla u\vert^2+\frac{1}{\epsilon}W(u)\Big)dx,\\
 &\mathscr{A}=\{u\in H^1(\Omega;\R^m):u=z,\;x\in\partial\Omega\},
 \end{split}
 \label{problem}
 \end{equation}
 where $\Omega\subset\R^n$, $n\geq 2$ is a bounded smooth domain and $z\in\R^m\setminus A$ is a fixed vector.

We let $u=u_\epsilon$ a minimizer of \eqref{problem}. We now introduce a special system of coordinates in a neighborhood of $\partial\Omega$ that we use for deriving upper and lower bounds for $J_\Omega^\epsilon(u_\epsilon)$.

For each $p\in\partial\Omega$ we let $\nu_p$ the unit exterior normal to $\partial\Omega$ at $p$. The smoothness of $\Omega$ implies that there is $h_0>0$ such that
\[x(p,h)=p-h\nu_p,\;\;p\in\partial\Omega,\;h\in(-h_0,h_0)\]
defines a diffeomorphism of $\partial\Omega\times(-h_0,h_0)$ onto $\{x\in\R^n:d(x,\partial\Omega)<h_0\}$. We have
\[\frac{\partial x(p,h)}{\partial(p,h)}=\begin{matrix}
1-hk_1& 0&\cdots &0&\\
0&1-hk_2& 0\cdots& 0&\\
\cdots\\
0&\cdots 0& 1-hk_{n-1}& 0& \\
0&\cdots &0 & 1&
\end{matrix}\]
and
\[\mathrm{det}(\frac{\partial x(p,h)}{\partial(p,h)})=\Pi_{j=1}^{n-1}(1-hk_j),\]
where $k_1,\ldots,k_{n-1}$ are the principal curvatures of $\partial\Omega$ at $p$.
We let $\sigma^+$, $\sigma^*$ and $\bar{u}_+$ as in $\S 2$. We have
\begin{lemma}$(\mathrm{Upper\;Bound})$
\label{upper}
There is a constant $C_0>0$ such that
\[J_\Omega^\epsilon(u)\leq\sigma^+\mathcal{H}^{n-1}(\partial\Omega)+C_0\epsilon,\]
for a minimizer $u$ of \eqref{problem}.
\end{lemma}
\begin{proof}
We define a comparison map $\tilde{u}:\bar{\Omega}\rightarrow\R^m$ by setting
\[\begin{split}
&\tilde{u}(x(p,h))=\bar{u}_+(\frac{h}{\epsilon}),\;\;p\in\partial\Omega,\;h\in[0,h_0),\\
&\tilde{u}(x)=\bar{u}_+(\frac{h_0}{\epsilon}),\;\;x\in\Omega,\;d(x,\partial\Omega)\geq h_0.
\end{split}\]
If $e_1,...,e_{n-1},\nu_p$ is an orthonormal basis with $e_1,...,e_{n-1}$ tangent to $\partial\Omega$ at $p$ in the direction of the principal curvatures, we have
\[\begin{split}
&\tilde{u}_{x_j}(x(p,h))=0,\;\;h\in(0,h_0),\;j=1,\ldots,n-1,\\
&\tilde{u}_{x_n}(x(p,h))=\frac{1}{\epsilon}\bar{u}_+^\prime(\frac{h}{\epsilon}),\;\;h\in(0,h_0).
\end{split}\]
It follows
\[\begin{split}
&J^\epsilon_{\{d(x,\partial\Omega)<h_0\}}(\tilde{u})=\frac{1}{\epsilon}\int_{\partial\Omega}\int_0^{h_0}
\Big(\frac{1}{2}\vert\bar{u}_+^\prime(\frac{h}{\epsilon})\vert^2
+W(\bar{u}_+(\frac{h}{\epsilon}))\Big)\vert\Pi_{j=1}^{n-1}(1-hk_j)\vert dhdp\\
&=\int_{\partial\Omega}\int_0^{\frac{h_0}{\epsilon}}\vert\bar{u}_+^\prime(s)\vert^2\vert\Pi_{j=1}^{n-1}(1-\epsilon sk_j)\vert dsdp\leq\sigma^+\mathcal{H}^{n-1}(\partial\Omega)+\frac{1}{2}C_0\epsilon,
\end{split}
\]
where $C_0>0$ is a constant and we have used $\frac{1}{2}\vert\bar{u}_+^\prime(s)\vert^2=W(\bar{u}_+(s))$.
We also have
\[\begin{split}
&J^\epsilon_{\{d(x,\partial\Omega)\geq h_0\}}(\tilde{u})=\frac{1}{\epsilon}W(\bar{u}_+(\frac{h_0}{\epsilon}))\int_{\{d(x,\partial\Omega)\geq h_0\}}dx\\
&\leq\frac{1}{2\epsilon}C_W^2\vert\bar{u}_+(\frac{h_0}{\epsilon})-a\vert^2\vert\Omega\vert
\leq\frac{1}{2\epsilon}C_W^2\bar{K}^2e^{-\frac{2\bar{k}h_0}{\epsilon}}\vert\Omega\vert\leq\frac{1}{2}C_0\epsilon,
\end{split}
\]
where we have used \eqref{W} and  $\vert\bar{u}_+(s)-a\vert\leq\bar{K}e^{-\bar{k}s}$
 for some $\bar{k}, \bar{K}>0$.
\end{proof}

Next we introduce some basic lemmas and the lower bound.

\begin{lemma}
\label{K-meas}
Let $u$ a minimizer of \eqref{problem} and $\delta>0$ a small number. Set $\Omega_{\delta}=\{x\in\Omega:\min_{a\in A}\vert u-a\vert\}\leq\delta\}$. Then there is a constant $C_1>0$ such that
\[\vert\Omega_{\delta}\vert\geq\vert\Omega\vert(1-\frac{C_1}{\delta^2}\epsilon),\]
where $\vert\vert$ stands for the $n$-dimensional Lebesgue measure.
\end{lemma}
\begin{proof}
Let $\Omega_\delta^c=\bar{\Omega}\setminus\Omega_\delta$. Then Lemma \ref{upper} and \eqref{outW} imply
\[\frac{1}{2}\frac{c_W^2\delta^2}{\epsilon}\vert\Omega_\delta^c\vert\leq\frac{1}{\epsilon}\int_\Omega W(u)dx\leq J_\Omega^\epsilon(u)
\leq\sigma^+\mathcal{H}^{n-1}(\partial\Omega)+C_0\epsilon.\]
It follows that there is $C_1>0$ such that
\[\begin{split}
&\vert\Omega_\delta^c\vert\leq\frac{C_1}{\delta^2}\vert\Omega\vert\epsilon,
\\
&\vert\Omega_\delta\vert\geq\vert\Omega\vert(1-\frac{C_1}{\delta^2}\epsilon).
\end{split}\]
\end{proof}

\begin{theorem}$(\mathrm{The\;Lower\;Bound})$
\label{lower}
Let $u$ a minimizer of \eqref{problem}. Then there is a constant $C>0$ such that
\begin{equation}
J_\Omega^\epsilon(u)\geq\sigma^+\mathcal{H}^{n-1}(\partial\Omega)(1-C\epsilon^\frac{1}{3}).
\label{lower1}
\end{equation}
\end{theorem}
\begin{proof}
For each $p\in\partial\Omega$ define $h_p$ by setting
\[\begin{split}
&h_p=h_0,\;\;\text{if}\;p-h\nu_p\in\Omega_\delta^c,\;h\in[0,h_0],\\
&h_p=\max\{h\in(0,h_0]:p-s\nu_p\in\Omega_\delta^c,\,s\in[0,h)\},\;\;\text{otherwise}.
\end{split}\]
Let $\alpha\in[0,1)$ a number to be chosen later and let
\[\mathscr{S}^\alpha=\{p\in\partial\Omega: h_p\geq\epsilon^\alpha\}.\]
We have
\[
\int_{\mathscr{S}^\alpha}\int_0^{\epsilon^\alpha}\vert\Pi_j^{n-1}(1-hk_j)\vert dhdp
\leq\int_{\mathscr{S}^\alpha}\int_0^{h_p}\vert \Pi_j^{n-1}(1-hk_j)\vert dhdp\leq\vert\Omega_\delta^c\vert
\leq\frac{C_1}{\delta^2}\vert\Omega\vert\epsilon.
\]
This and $\vert\Pi_j^{n-1}(1-hk_j)\vert\geq 1-C_k\epsilon^\alpha$, ($C_k>0$ a constant that depends on the curvatures of $\partial\Omega$) imply the existence of $C_2>0$ such that
\begin{equation}
\begin{split}
&\mathcal{H}^{n-1}(\mathscr{S}^\alpha)\epsilon^\alpha(1-C_k\epsilon^\alpha)\leq\frac{C_1}{\delta^2}\vert\Omega\vert\epsilon ,\\
&\Rightarrow \mathcal{H}^{n-1}(\mathscr{S}^\alpha)\leq\frac{C_2}{\delta^2}\epsilon^{1-\alpha}\mathcal{H}^{n-1}(\partial\Omega).
\end{split}
\label{S}
\end{equation}
On the other hand Lemma \ref{lower-sigma}, Lemma \ref{basic} and \eqref{S} imply
\[\begin{split}
&J_\Omega^\epsilon(u)\geq(\sigma^+-\frac{1}{2}C_W\delta^2)\int_{\partial\Omega\setminus\mathscr{S}^\alpha}\min\{1,
\vert\Pi_j^{n-1}(1-\epsilon^\alpha k_j)\vert\}dp\\
&\geq(\sigma^+-\frac{1}{2}C_W\delta^2)\int_{\partial\Omega\setminus\mathscr{S}^\alpha}( 1-C_k\epsilon^\alpha)dp\\
&\geq(\sigma^+-\frac{1}{2}C_W\delta^2)\mathcal{H}^{n-1}(\partial\Omega)(1-\frac{C_2}{\delta^2}\epsilon^{1-\alpha})( 1-C_k\epsilon^\alpha).
\end{split}\]
 It follows
\[J_\Omega^\epsilon(u)\geq\sigma^+\mathcal{H}^{n-1}(\partial\Omega)(1-C_3(\delta^2+\epsilon^\alpha+\frac{\epsilon^{1-\alpha}}{\delta^2})),\]
for some constant $C_3>0$. This implies \eqref{lower1} for $\alpha=\frac{1}{3}$, $\delta^2=\epsilon^\frac{1}{3}$.
\end{proof}
Let $\delta_0\in(0,\min_{i\neq j}\frac{1}{2}\vert a_i-a_j\vert)$ and $\delta\in(0,\delta_0)$ be fixed.
A consequence of the upper bound, lower bound and the density estimate ( see for example Theorem 5.2 in\cite{afs}) is that a minimizer $u$ remains near just one of the zeros of $W$ throughout $\Omega$. We state now the main result of this section
\begin{theorem}
\label{nearone} Let $u$, $(u=u_\epsilon)$ a minimizer of problem \eqref{problem}.
There exist $a_+\in A$ and positive constants $k, K$ and $C$ such that
\[
\vert u(x)-a_+\vert\leq Ke^{-\frac{k}{\epsilon}(d(x,\partial\Omega)-C\epsilon^\frac{1}{3(n-1)})},\;\;x\in\Omega.
\]
\end{theorem}
\begin{proof}
1. Assume $x_0\in\Omega_\delta^c$. Then the gradient bound
\[\vert\nabla u\vert\leq \frac{M}{\epsilon}\]
implies
\begin{equation}
\min_j\vert u(x)-a_j\vert>\frac{\delta}{2},\;\;x\in B_{\frac{\delta\epsilon}{2M}}(x_0).
\label{deltamezzi}
\end{equation}
Let $v:\frac{\Omega}{\epsilon}\rightarrow\R^m$ be defined by $v(\frac{x}{\epsilon})=u(x)$. Then the minimality of $u$ implies that $v$ minimizes $\int_{\frac{\Omega}{\epsilon}}(\frac{1}{2}\vert\nabla v\vert^2+W(v))dy$ and \eqref{deltamezzi} implies that $v$ satisfies
\[\min_j\vert v(y)-a_j\vert>\frac{\delta}{2},\;\;y\in B_{\frac{\delta}{2M}}(y_0).\]
Then arguing as in the proof of Lemma 5.5 in \cite{afs} we obtain, via the density estimate, with $c^\prime(\delta)$ a constant,
\[\vert B_{r}(y_0)\cap\{y:\min_j\vert v(y)-a_j\vert>\delta\}\vert\geq(\delta_0-\delta)c^\prime(\delta)r^{n-1},\;\; r\leq d(y_0,\partial(\frac{\Omega}{\epsilon})),\]
which, in terms of $u$, becomes
\[\epsilon^{-n}\vert B_{\epsilon r}(x_0)\cap\Omega_\delta^c\vert\geq(\delta_0-\delta)c^\prime(\delta)r^{n-1},\;\;\epsilon r\leq d(x_0,\partial\Omega).\]
It follows
\begin{equation}
J^\epsilon_{B_{\epsilon r}(x_0)}(u)\geq w_\delta(\delta_0-\delta)c^\prime(\delta)(\epsilon r)^{n-1}=C_\delta(\epsilon r)^{n-1},\;\;\epsilon r\leq d(x_0,\partial\Omega),
\label{extra}
\end{equation}
where $w_\delta=\min_{x\in\Omega_\delta^c}W(u(x))>0$ and $C_\delta=w_\delta(\delta_0-\delta)c^\prime(\delta)$.

2. From the proof of Theorem \ref{lower} we see that the right end side of \eqref{lower1} is an estimate of the energy contained in an $\epsilon^\frac{1}{3}-$neighborhood of $\partial\Omega$. Therefore we can add the energy
$J_{B_{\epsilon r}(x_0)}$ and improve the lower bound given by \eqref{lower1} provided $x_0$ and $r$ satisfy the
condition
\[d(x_0,\partial\Omega)>\epsilon^\frac{1}{3}+\epsilon r.\]
In order to reach a contradiction with the upper bound given by Lemma \ref{upper} we choose $r$ by imposing the condition
\[C_\delta(\epsilon r)^{n-1}=2\sigma^+\mathcal{H}^{n-1}(\partial\Omega) C\epsilon^\frac{1}{3}\;\Leftrightarrow\; r=(\frac{2\sigma^+\mathcal{H}^{n-1}(\partial\Omega) C}{C_\delta})^\frac{1}{n-1}\epsilon^{-(1-\frac{1}{3(n-1)})},\]
where $C>0$ is the constant in \eqref{lower1}. With this choice of $r$ we obtain that the existence of a point
$x_0\in\Omega_\delta^c$ that satisfies $d(x_0,\partial\Omega)>\epsilon^\frac{1}{3}+C^\prime\epsilon^{\frac{1}{3(n-1)}}$ ($C^\prime=(\frac{2\sigma^+\mathcal{H}^{n-1}(\partial\Omega) C}{C_\delta})^\frac{1}{n-1}$)
implies
\[J_\Omega^\epsilon(u)\geq\sigma^+\partial\Omega(1+C\epsilon^\frac{1}{3})\]
that collides with \ref{upper}. This and the smoothness of the minimizer $u$ imply that there is $a_0\in A$ such that
\begin{equation}
d(x_0,\partial\Omega)\geq(2C^\prime+1)\epsilon^{\frac{1}{3(n-1)}}\;\Rightarrow\;\vert u(x)-a_0\vert<\delta.
\label{near a0}
\end{equation}

3. $a_0\in\{a\in A: a=a_+\}$. If this is not the case, from Lemmas \ref{lower-sigma},  \ref{basic} and \ref{upper} we have 
\[(\sigma^*-\frac{1}{2}C_W\delta^2)\mathcal{H}^{n-1}(\partial\Omega)(1-C_K\epsilon^{\frac{1}{3(n-1)}})
\leq\sigma^+\mathcal{H}^{n-1}(\partial\Omega)+C_0\epsilon,\]
where $C_K>0$ is a constant that depends on $C^\prime$ and from the curvature of $\partial\Omega$.
Since $\sigma^*>\sigma^+$ this inequality cannot hold if $\delta$ and $\epsilon$ are sufficiently small. This proves the claim.

From \eqref{near a0} with $a_0=a_+$ and linear elliptic theory the exponential inequality follows. The proof is complete.
\end{proof}

\section{ Example 2}\label{ex2}

We consider the minimization problem
 \begin{equation}
 \begin{split}
 &\min_{u\in\mathscr{A}} J_\Omega^\epsilon(u),\;\;J_\Omega^\epsilon(u)=\int_\Omega\Big(\frac{\epsilon}{2}\vert\nabla u\vert^2+\frac{1}{\epsilon}W(u)\Big)dx,\\
 &\mathscr{A}=\{u\in H^1(\Omega;\R^m):u=g_\epsilon,\;x\in\partial\Omega\}.
 \end{split}
 \label{problem2}
 \end{equation}
 under specific assumptions on the potential $W:\R^m\rightarrow\R$, $m\geq 1$, the domain $\Omega\subset\R^2$ and the Dirichlet data $g_\epsilon$ that we now describe
 \begin{description}
 \item[h$_1$] $W$ is a $C^2$ nonnegative function that vanishes on a finite set: $\{W=0\}=A$,

  $A=\{a_1,\ldots,a_N\}$, $N\geq 2$.

     Moreover
 \[\begin{split}
 &\liminf_{\vert u\vert\rightarrow+\infty} W(u)>0,\\
 &W_u(u)\cdot u>0,\quad\text{for}\;\vert u\vert>M,\;\text{some}\;M>0,
 \end{split}\]
 and
 \[\xi^\top W_{uu}(a_j)\xi\geq 2c^2\vert\xi\vert^2,\;j=1,\ldots,N,\;\text{some}\;c>0.\]
 \end{description}

 From Lemma \ref{E-conn} we have that $\mathbf{h}_1$ implies that given $a_-\in A$ there exists $a_+\in A\setminus\{a_-\}$ and a map $\bar{u}:\R\rightarrow\R^m$, connecting $a_\pm$ at $\pm\infty$, that minimizes
 $J_\R(v)=\int_\R(\frac{1}{2}\vert v^\prime\vert^2+W(v))ds$ in the class of $H_{\mathrm{loc}}^1$ maps that connect $a_-$ to $A\setminus\{a_-\}$ (see Theorem 2.1 in \cite{afs}). In general $a_+$ is not uniquely determined and we let $A_+\subset A$ the set of the elements in $A$ that can be identified with $a_+$.  There exists $\eta>0$ such that, for $a\in A\setminus A_+$ and $\delta>0$ sufficiently small, for each map $v:(s_1,s_2)\rightarrow\R^m$, $-\infty\leq s_1<s_2\leq+\infty$, that satisfies
 \[\lim_{s\rightarrow s_1}v(s)\in B_\delta(a_-),\quad\;\lim_{s\rightarrow s_2}v(s)\in B_\delta(a),\]
 it results
 \[J_{(s_1,s_2)}(v)\geq\sigma+\eta,\]
 where $\sigma=\int_\R\vert\bar{u}^\prime\vert^2ds$ is
 the energy of $\bar{u}$.

Concerning $\Omega$ and the Dirichlet data $g_\epsilon$ we assume
\begin{description}
\item[h$_2$] $\Omega\subset\R^2$ is open, bounded, connected and of class $C^{1,\alpha}$ and satisfies the condition
\begin{equation}
\Omega\cap[0,l]\times[0,h]=[0,l]\times(0,h),
\label{Omega}
\end{equation}
where $l,h$ are positive constants.

See Figure \ref{OOmega} for an example of a set $\Omega$ that satisfies $\mathbf{h}_2$.

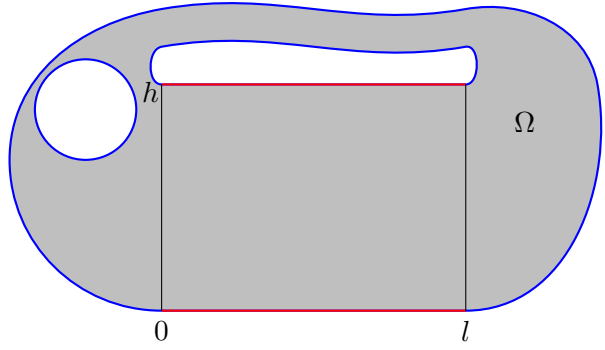
\begin{figure}
\begin{center}
\begin{tikzpicture}


\path [fill=lightgray] (0,0) -- (4,0) to [out=0,in=280] (5.732,3) to [out=100,in=10] (4,4)  to [out=190,in=10] (0,4)
to [out=190,in=90] (-2,2) to [out=270,in=180] (0,0);

\path [fill=white] (0,3) -- (4,3) to [out=0,in=10] (4,3.5) to [out=190,in=10] (0,3.5) to [out=190,in=180] (0,3);
\path [fill=white]  (-1,2.6666) circle [radius=0.6666];

\draw [blue, thick] (0,0) -- (4,0) to [out=0,in=280] (5.732,3) to [out=100,in=10] (4,4)  to [out=190,in=10] (0,4)
to [out=190,in=90] (-2,2) to [out=270,in=180] (0,0);

\draw [blue, thick] (0,3) -- (4,3) to [out=0,in=10] (4,3.5) to [out=190,in=10] (0,3.5) to [out=190,in=180] (0,3);
\draw [blue, thick]  (-1,2.6666) circle [radius=0.6666];

\draw [red, thick] (0,0) -- (4,0);
\draw [red, thick] (0,3) -- (4,3);

\draw [thin](0,0)--(0,3);
\draw [thin](4,0)--(4,3);

\node[below] at (0,0) {$0$};
\node[below] at (4,0) {$l$};
\node[left] at (.1,2.9) {$h$};
\node[right] at (4.5,2.5) {$\Omega$};
\end{tikzpicture}
\end{center}
\caption{An example of $\Omega$  case $l>h$.}
\label{OOmega}
\end{figure}
\item[h$_3$] $g_\epsilon$ satisfies $\vert g_\epsilon\vert\leq M$ and
\begin{equation}
\begin{split}
&g_\epsilon(x,0)=g_\epsilon(x,h)=a_+,\;\;x\in(C_0\epsilon,l-C_0\epsilon),\;\;\text{some}\;C_0>0,\\
&g_\epsilon(x,y)=a_-,\;\;(x,y)\in\partial\Omega\setminus(0,l)\times\{0,h\}=:\partial\Omega^-,
\\
&\vert g_{\epsilon,x}(x,y)\vert\leq \frac{C}{\epsilon},\;\;x\in(0,C_0\epsilon)\cup(l-C_0\epsilon,l)\;\;\vert g_{\epsilon,x}(\cdot,y)\vert_\alpha\leq\frac{C}{\epsilon^{1+\alpha}},\;y=0,h.
\end{split}
\label{g-def}
\end{equation}
\end{description}

\begin{lemma}
\label{Csueps}
Minimizers $u_\epsilon$ of problem \eqref{problem2} (actually $H^1(\Omega;\R^m)$ critical points ) satisfy the estimates
\begin{equation}
\begin{split}
&\Vert u_\epsilon\Vert_{L^\infty}<M,\\
&\Vert\nabla u_\epsilon\Vert_{L^\infty}<\frac{C^{\prime\prime}}{\epsilon}.
\end{split}
\label{bounds}
\end{equation}
\end{lemma}
\begin{proof}
By linear elliptic theory $u_\epsilon\in C^2(\Omega;\R^m)$. Set $v_\epsilon=\vert u_\epsilon\vert^2$. Then
\[\Delta v_\epsilon=2W_u(u_\epsilon)\cdot u_\epsilon+2\vert\nabla u_\epsilon\vert^2>0,\;\;\text{for}\;\vert u_\epsilon\vert>M.\]
Hence $\max\vert u_\epsilon\vert^2\leq M^2$ if $v_\epsilon$ attains its max in the interior of $\Omega$.
\noindent
On the other hand by $\mathbf{h}_3$ $\max_{\partial\Omega}\vert u_\epsilon\vert\leq M$ and \eqref{bounds}$_1$ follows.
\noindent
For the gradient bound rescale the solution of the Euler-Lagrange equation
\[\epsilon\Delta u_\epsilon=\frac{1}{\epsilon}W_u(u_\epsilon),\]
by $\frac{z}{\epsilon}$ and denote by $\tilde{u}, \tilde{g}$ the rescaled $u_\epsilon,
g_\eps$. Then by linear theory (see (8.87) in \cite{GT})
\[\vert\tilde{u}\vert_{1+\alpha}\leq C^\prime(\Vert\tilde{u}\Vert_{L^\infty}+\vert\tilde{g}\vert_{1,\alpha})\leq C^{\prime\prime}.\]
Hence
\[\Vert\nabla \tilde{u}\Vert_{L^\infty}\leq C^{\prime\prime},\]
and therefore \eqref{bounds}$_2$ follows.
\end{proof}

\subsection{$\frac{h}{l}>1$ The boundary layer case}
In this section we analyze in detail the structure of the minimizers of $J_\Omega(u)$ under the assumption
\begin{equation}
l<h.
\label{l<h}
\end{equation}

 We will show that, for $\epsilon>0$ small, minimizers of \eqref{problem2} remain near $a_-$ in the whole $\Omega$ aside from a boundary layer near the upper and lower boundary of the rectangle $(0,l)\times(0,h)$ where the condition $u=a_+$ is imposed. We will also estimate from below the thickness of the boundary layer (cfr. Theorem \ref{thick}).

 Set
 \[\partial^+\Omega=(0,l)\times\{0,h\}.\]

 Then we have
 \begin{theorem}
\label{BL-exist}
There are positive constants $k,C$ and $C^1$ independent of $\epsilon$ and such that
\begin{equation}
\vert u(z)-a_-\vert\leq
C e^{-\frac{k}{\epsilon}(d(z,\partial^+\Omega)-C^1\epsilon^\frac{1}{4}\vert\ln{\epsilon}\vert^\frac{3}{2})^+},\;\;z\in\Omega.
\label{exp-bound}
\end{equation}
\end{theorem}
The proof of Theorem \ref{BL-exist} is based on the tight lower$/$upper energy bounds and on the vector analog of the Caffarelli-Cordoba density estimate (see [\cite{afs} Theorem 5.2]).
\vskip.2cm
Set $J_\R^\epsilon(v)=\int_\R(\frac{\epsilon}{2}\vert v^\prime\vert^2+\frac{1}{\epsilon}W(v))ds$ and note that
\[\sigma= J_\R^\epsilon(\bar{u}(\frac{\cdot}{\epsilon}))=J_\R(\bar{u}).\]
\subsubsection{The Upper Bound}
\begin{proposition}
\label{UB2}
There is a constant $C_1>0$ such that, if
 $u:\Omega\rightarrow\R^m$ is a minimizer of the problem \eqref{problem2} then
\begin{equation}
J_\Omega^\epsilon(u)\leq 2l\sigma+C_1\epsilon\vert\ln{\epsilon}\vert^3.
\label{UB*}
\end{equation}
\end{proposition}
\begin{remark}
We don't believe that the logarithm  can be removed.
\end{remark}
\begin{proof}
We show that there exist a test map $\tilde{u}:\Omega\rightarrow\R^m$ that satisfies the bound \eqref{UB*}.
Let $\eta>0$ a small number to be chosen later and set
\[\begin{split}
&L_1=[C_0\epsilon,l-C_0\epsilon]\times[0,2\eta],\;L_2=[C_0\epsilon,l-C_0\epsilon]\times[h-2\eta,h],\\
&T^1=\{(x,y):0\leq y\leq\frac{2\eta}{C_0\epsilon}x,\;x\in[0,C_0\epsilon]\}.
\end{split}\]
We denote by $T^2$ the reflection of $T^1$ in the $x$ axis, by $T^3$ the reflection of $T^1$ in the $y$ axis, by $T_4$ the reflection of $T^2$ in the $y$ axis and by $T_\tau^j$, $T^j$ subjected to the translation $\tau\in\R^2$.
Set
\[B=L_1\cup L_2\cup T^1\cup T_{(0,h)}^2\cup T_{(l,0)}^3\cup T_{(l,h)}^4.\]

 We set $\tilde{u}(x,y)=a_-$ for $(x,y)\in\Omega\setminus B$. To define $\tilde{u}$ in $L_1\cup L_2$ we use the map $\bar{v}:[-\eta,\eta]\rightarrow\R^m$ given by
\[\bar{v}(s)=\left\{\begin{array}{l}
\frac{1}{C_2\epsilon}\Big((-\eta+C_2\epsilon-s)a_- +(\eta+s)\bar{u}(C_2-\frac{\eta}{\eps})\Big),\;\;-\eta\leq s\leq -\eta+C_2\epsilon,\\
\bar{u}(\frac{s}{\epsilon}),\;\;-\eta+C_2\epsilon\leq s\leq\eta-C_2\epsilon,\\
\frac{1}{C_2\epsilon}\Big((\eta-s)\bar{u}(\frac{\eta}{\eps}-C_2)+(s-\eta+C_2\epsilon)a_+)\Big),
\;\;\eta-C_2\epsilon\leq s\leq\eta.
\end{array}
\right.\]
From ${\mathbf{h}}_1$ it follows that there are positive constants $k$ and $K$ such that
\begin{equation}
\begin{split}
&\vert\bar{u}(s)-a_-\vert\leq Ke^{ks},\;\;s\leq 0,\\
&\vert\bar{u}(s)-a_+\vert\leq Ke^{-ks},\;\;s\geq 0.
\end{split}
\label{kK}
\end{equation}
We take $\eta=\frac{1}{2k}\epsilon\vert\ln{\epsilon}\vert$. Then \eqref{kK} implies
\[\begin{split}
&\vert\bar{u}(C_2-\frac{\eta}{\eps})-a_-\vert\leq Ke^{C_2}e^{\ln{\epsilon^\frac{1}{2}}}=C_3\epsilon^\frac{1}{2},\\
&\vert\bar{u}(\frac{\eta}{\eps}-C_2)-a_+\vert\leq Ke^{C_2}e^{\ln{\epsilon^\frac{1}{2}}}=C_3\epsilon^\frac{1}{2}.
\end{split}\]
This, \eqref{W}, the definition of $\bar{v}$ and a standard computation yield
\begin{equation}
J^\epsilon_{(-\eta,-\eta+C_2\epsilon)}(\bar{v}), J^\epsilon_{(\eta-C_2\epsilon,\eta)}(\bar{v})
\leq\frac{1}{2}(\frac{C_3^2}{C_2}+C_W^2C_3^2C_2)\epsilon=C_4\epsilon.
\label{jv-est}
\end{equation}
This and the obvious inequality $J^\epsilon_{(-\eta+C_2\epsilon,\eta-C_2\epsilon)}(\bar{v})\leq\sigma$ yield
\begin{equation}
J^\epsilon_{(-\eta,\eta)}(\bar{v})\leq\sigma+2C_4\epsilon.
\label{jbarv}
\end{equation}
We set
\begin{equation}
\tilde{u}(x,y)=\left\{
\begin{array}{l}
\bar{v}(y-h+\eta),\;(x,y)\in L_2,\\
\bar{v}(\eta-y),\;(x,y)\in L_1.
\end{array}
\right.
\label{utildeL}
\end{equation}
It remains to define $\tilde{u}$ on $T^1\cup T_{(0,h)}^2\cup T_{(l,0)}^3\cup T_{(l,h)}^4$. We define $\tilde{u}$
 on $T^1$. The definition on $T_{(0,h)}^2, T_{(l,0)}^3,T_{(l,h)}^4$ is similar. Let $(r,s)\in\R^2$ coordinates with respect to a positively oriented system with origin in $(0,0)$ and the $s$ axis coinciding with the line $y=\frac{2\eta}{C_0\epsilon}x$. We define $\tilde{u}\vert_{T^1}$ on each line parallel to the $s$ axis by linear interpolation between the values on the intersection of the line with the boundary of $T^1$. Set $\lambda=(C_0^2\epsilon^2+4\eta^2)^\frac{1}{2}$, $\alpha=\frac{C_0\epsilon}{\lambda}$ and $\beta=\frac{2\eta}{\lambda}$. With this definitions it results
 \begin{equation}
 \tilde{u}(r,s)=\frac{\lambda-\frac{\beta}{\alpha}r-s}{\lambda-\frac{r}{\alpha\beta}}g(\frac{r}{\beta},0)
 +\frac{s-\frac{\alpha}{\beta}r}{\lambda-\frac{r}{\alpha\beta}}\bar{v}(-\eta+\frac{r}{\alpha}),\;\;
 s\in[\frac{\alpha}{\beta}r,\lambda-\frac{\beta}{\alpha}r],\;r\in[0,\alpha\beta\lambda].
 \label{Tdef}
 \end{equation}
 Since $\tilde{u}$ is uniformly bounded for $\epsilon\rightarrow 0^+$ the same is true for $W(\tilde{u})$. It follows
 \begin{equation}
  \frac{1}{\epsilon}\int_{T^1}W(\tilde{u})\leq CC_0\eta=C\epsilon\vert\ln{\epsilon}\vert.
  \label{Tpot-est}
 \end{equation}
From \eqref{Tdef} we have
\[\tilde{u}_s(r,s)=\frac{1}{\lambda-\frac{r}{\alpha\beta}}(\bar{v}(-\eta+\frac{r}{\alpha})-g(\frac{r}{\beta},0)),\]
and therefore
\[\epsilon\int_{T^1}\vert\tilde{u}_s\vert^2
=\epsilon\int_0^{\alpha\beta\lambda}
\int_{\frac{\alpha}{\beta}r}^{\lambda-\frac{\beta}{\alpha}r}\vert\tilde{u}_s\vert^2dsdr
=\epsilon\int_0^{\alpha\beta\lambda}\frac{\alpha\beta}{\alpha\beta\lambda-r}
\vert\bar{v}(-\eta+\frac{r}{\alpha})-g(\frac{r}{\beta},0)\vert^2dr.\]
Observe that from $\beta\lambda=2\eta$ and $\alpha\lambda=C_0\epsilon$ it follows
\[r=\alpha\beta\lambda\;\;\Rightarrow\;\bar{v}(-\eta+\frac{r}{\alpha})=\bar{v}(\eta)
=g(\frac{r}{\beta},0))=g(C_0\epsilon,0)=a_+.\]
This and the fact that $\vert\bar{v}^\prime\vert,\vert g_x\vert\leq\frac{C}{\epsilon}$ imply
\begin{equation}
\begin{split}
&\vert\bar{v}(-\eta+\frac{r}{\alpha})-g(\frac{r}{\beta},0)\vert
\leq\vert\bar{v}(-\eta+\frac{r}{\alpha})-a_+\vert+\vert g(\frac{r}{\beta},0)-a_+\vert\\
&\leq\frac{C}{\epsilon}(\vert 2\eta-\frac{r}{\alpha}\vert+\vert\frac{r}{\beta}-\alpha\lambda\vert)
=\frac{C}{\epsilon}\frac{\alpha+\beta}{\alpha\beta}(\alpha\beta\lambda-r).
\end{split}
\label{gx-vprime}
\end{equation}
With this estimate we obtain
\[\begin{split}
&\epsilon\int_{T^1}\vert\tilde{u}_s\vert^2\leq\epsilon\int_0^{\alpha\beta\lambda}
\frac{C^2}{\epsilon^2}\frac{(\alpha+\beta)^2}{\alpha\beta}(\alpha\beta\lambda-r) dr\\
&=\frac{C^2}{2\epsilon}\alpha\beta\lambda^2(\alpha+\beta)^2\leq\frac{C^2}{\epsilon}\alpha\beta\lambda^2\leq C\epsilon\vert\ln{\epsilon}\vert.
\end{split}\]
After some manipulation we obtain
\[\begin{split}
&\tilde{u}_r(r,s)=\frac{1}{\alpha\beta}(\alpha^2\lambda-s)\frac{g(\frac{r}{\beta},0)-\bar{v}(-\eta+\frac{r}{\alpha})}
{(\lambda-\frac{r}{\alpha\beta})^2}\\
&+\frac{1}{\beta}\frac{\lambda-\frac{\beta}{\alpha}r-s}{\lambda-\frac{r}{\alpha\beta}}g_x(\frac{r}{\beta},0)
 +\frac{1}{\alpha}\frac{s-\frac{\alpha}{\beta}r}{\lambda-\frac{r}{\alpha\beta}}\bar{v}^\prime(-\eta+\frac{r}{\alpha}).
\end{split}\]
This and \eqref{gx-vprime} imply
\[\vert\tilde{u}_r\vert^2\leq\frac{C}{\epsilon^2}\frac{(\alpha^2\lambda-s)^2
+\alpha^2(\lambda-\frac{\beta}{\alpha}r-s)^2+\beta^2(s-\frac{\alpha}{\beta}r)^2}{(\alpha\beta\lambda-r)^2}.\]
Observing that
\[\int_{\frac{\alpha}{\beta}r}^{\lambda-\frac{\beta}{\alpha}r}
\Big((\alpha^2\lambda-s)^2
+\alpha^2(\lambda-\frac{\beta}{\alpha}r-s)^2+\beta^2(s-\frac{\alpha}{\beta}r)^2\Big)ds
=\frac{1+\alpha^6+\beta^6}{3\alpha^3\beta^3}(\alpha\beta\lambda-r)^3,\]
we finally obtain
\[\begin{split}
&\epsilon\int_{T^1}\vert\tilde{u}_r\vert^2
=\epsilon\int_0^{\alpha\beta\lambda}
\int_{\frac{\alpha}{\beta}r}^{\lambda-\frac{\beta}{\alpha}r}\vert\tilde{u}_r\vert^2dsdr\\
&\leq\frac{C}{\epsilon\alpha^3\beta^3}\int_0^{\alpha\beta\lambda}(\alpha\beta\lambda-r)dr
=\frac{C\lambda^2}{2\epsilon\alpha\beta}\leq C\epsilon\vert\ln{\epsilon}\vert^3.
\end{split}\]
The proof is complete.
\end{proof}
\subsubsection{The Lower Bound.}
We now derive a lower bound for the energy of a minimizer $u:\Omega\rightarrow\R^m$ of problem \eqref{problem2}.

\begin{proposition}
\label{LB*}
There exist $\epsilon_0>0$ and  a constant $\tilde{C}>0$ independent of $\epsilon\in(0,\epsilon_0]$ such that, if $u:\Omega\rightarrow\R^m$ is a minimizer of problem \eqref{problem2}, then
\begin{equation}
J_\Omega^\epsilon(u)\geq\int_0^l\int_0^h\Big(\frac{\epsilon}{2}\vert\frac{\partial u}{\partial y}\vert^2+\frac{1}{\epsilon}W(u)\Big)dydx\geq 2l\sigma-\tilde{C}\epsilon^\frac{1}{2}.
\label{LB*1}
\end{equation}
\end{proposition}
\begin{remark}
The following lower bound is derived below (see \eqref{LB-ref})
\[J_\Omega^\epsilon(u)\geq\int_0^l\int_0^h\Big(\frac{\epsilon}{2}\vert\frac{\partial u}{\partial y}\vert^2+\frac{1}{\epsilon}W(u)\Big)dydx \geq 2\sigma l-\tilde{C}^\prime\epsilon,\]
which is an improvement over \eqref{LB*1}. Notice that only half of the gradient is involved.
\end{remark}
\begin{proof}
1. Given $y\in(0,h)$, let $\Sigma_y=\{(x^\prime,y^\prime)\in\Omega: y^\prime=y\}$ and let $Y\subset(0,h)$ be the set
\[Y=\{y\in(0,h):\Sigma_y\cap\{\min_{a\neq a_-}\vert u-a\vert\leq\delta\}=\emptyset\},\]
where $\delta=\delta_\epsilon$ is to be chosen later. That is the set of sections $\Sigma_y$ on which $u$ is $\delta$ away from all $a\neq a_-$.
To estimate the measure $\mathcal{H}^1(Y)$ from below we note that $y\in(0,h)\setminus Y$ implies the existence of $a\neq a_-$ and  $(x(y),y)$ such that $\vert u(x(y),y)-a\vert<\delta$. From this, \eqref{UB*} and Lemma \ref{lower-sigma} we obtain
\[\begin{split}
&2l\sigma+C_1\epsilon\vert\ln{\epsilon}\vert^3
\geq
\int_{(0,h)\setminus Y}\int_{\Sigma_y}(\frac{\epsilon}{2}\vert u_x(x,y)\vert^2+\frac{1}{\epsilon}W(u(x,y)))dx dy\\
&\geq(2\sigma-C_W\delta^2)(h-\mathcal{H}^1(Y)).
\end{split}\]
It follows, for $\delta>0$ small
\[\mathcal{H}^1(Y)\geq h-l-\mathrm{O}(\delta^2)
\geq\frac{1}{2}(h-l).\]

2. Among the sections $\Sigma_y$ considered in the previous step (a large set) consider those that contain a subset on which $u$ is also $\delta$-away from $a_-$. We will show that this subset has small measure for most of these sections.
Set $\Sigma_y^*=\{(x,y)\in\Sigma_y:\vert u(x,y)-a_-\vert>\delta\}$. Then the definition of $Y$ implies
\[\min_{a\in A}\vert u-a\vert>\delta,\;\;(x,y)\in\Sigma_y^*,\;y\in Y,\]
and therefore, from Lemma \ref{Wlemma} and Proposition \ref{UB2} we get
\begin{equation}
\begin{split}
&\frac{1}{2\epsilon}C_W^2\delta^2\int_Y\mathcal{H}^1(\Sigma_y^*)\leq 2l\sigma +C_1\epsilon\vert\ln{\epsilon}\vert^3,\\
&\Rightarrow \int_Y\mathcal{H}^1(\Sigma_y^*)\leq C^*\frac{\epsilon}{\delta^2},\;\;C^*=8\frac{l\sigma}{C_W^2}.
\end{split}
\label{bo}
\end{equation}
Let $Y^*=\{y\in Y:\mathcal{H}^1(\Sigma_y^*)\geq \frac{2}{\mathcal{H}^1(Y)}C^*\frac{\epsilon}{\delta^2}\}$. Then \eqref{bo} yields
$\mathcal{H}^1(Y^*)\leq\frac{1}{2}\mathcal{H}^1(Y)$,
and in turn by step 1.
\[\mathcal{H}^1(Y\setminus Y^*)\geq\frac{1}{4}(h-l).\]
For each $y\in Y\setminus Y^*$ we have that $\mathcal{H}^1(\Sigma_y^*)< \frac{2}{\mathcal{H}^1( Y)}C^*\frac{\epsilon}{\delta^2}$. It follows
\begin{equation}
y\in Y\setminus Y^*\;\;\Rightarrow\;\mathcal{H}^1(\{x\in(0,l):\vert u(x,y)-a_-\vert\leq\delta\})
>l-\frac{2}{\mathcal{H}^1(Y}C^*\frac{\epsilon}{\delta^2}.
\label{good}
\end{equation}

3. From step 2. there exist $y_1,y_2\in(0,h)$, $y_2-y_1\geq\frac{1}{4}(h-l)$, that satisfy \eqref{good}. The boundary condition implies $u(x,h)=a_+$ for $x\in(C_0\epsilon,l-C_0\epsilon)$. Hence from \eqref{good} we have
\[\mathcal{H}^1(\{ x\in(C_0\epsilon,l-C_0\epsilon):\vert u(x,y_2)-a_-\vert\leq\delta\})>l-2C_0\epsilon-\frac{2}{\mathcal{H}^1(Y)}C^*\frac{\epsilon}{\delta^2}.\]
This and Lemma \ref{lower-sigma} imply
\begin{equation}
\begin{split}
&J^\epsilon_{(C_0\epsilon,l-C_0\epsilon)\times(y_2,h)}(u)\geq\int_{\{ x\in(C_0\epsilon,l-C_0\epsilon):\vert u(x,y_2)-a_-\vert\leq\delta\}}J^\epsilon_{(y_2,h)}(u(x,\cdot))dx\\
&\geq(\sigma-\frac{1}{2}C_W\delta^2)(l-2C_0\epsilon-\frac{2}{\mathcal{H}^1(Y)}C^*\frac{\epsilon}{\delta^2}).
\end{split}
\label{halph}
\end{equation}
Since a similar estimate holds for $J^\epsilon_{(C_0\epsilon,l-C_0\epsilon)\times(0,y_1)}(u)$, the lower bound \eqref{LB*1} follows from \eqref{halph} with $\delta=\epsilon^\frac{1}{4}$ and from the fact that $J^\epsilon_{(0,y_1)}(u(x,\cdot))$ and $J^\epsilon_{(y_2,h)}(u(x,\cdot))$ only account for the derivative $u_y$.
The proof is complete.
\end{proof}
\begin{remark}
\label{E-in-D}
The lower bound in Proposition \ref{LB*} is an estimate from below of the energy of the minimizer $u$ in the set $D=[0,l]\times([0,y_1]\cup[y_2,h])$ which is quite tight with respect the upper bound \eqref{UB*} on the whole domain. Hence we expect that, in $\Omega\setminus D$, the minimizer to be close to $a_-$ uniformly.
\end{remark}
 Set $\delta_0=\frac{1}{2}\min_{i\neq j}\vert a_i-a_j\vert$.
\begin{lemma}
\label{nearone*}
Fix $\delta\in(0,\delta_0]$ then there is a constant $C_\delta>0$ such that
\begin{equation}
z\in\Omega,\;d(z,D))\geq C_\delta\epsilon^\frac{1}{2}\;\;\Rightarrow\;
\vert u(z)-a_-\vert\leq\delta.
\label{uptoB}
\end{equation}
Moreover there is a constant $k>0$ such that
\begin{equation}
\vert u(z)-a_-\vert\leq \delta e^{-\frac{k}{\epsilon}(d(z,D)-C_\delta\epsilon^\frac{1}{2})^+},\;\;z\in\Omega\setminus D.
\label{exp-refin}
\end{equation}
\end{lemma}
\begin{proof}
1. We first establish \eqref{uptoB}. Suppose that $\min_{a\in A}\vert u(z_0)-a\vert>\delta$ for some $z_0\in\Omega\setminus D$. From the gradient bound $\vert\nabla u\vert\leq\frac{M}{\epsilon}$ it follows that $\min_{a\in A}\vert u(z)-a\vert>\frac{\delta}{2}$ for $z\in B_{\frac{\delta}{2}}(z_0)$. Then by the Caffarelli-Cordoba (vector version) density estimate we obtain
\begin{equation}
J^\epsilon_{B_{\epsilon r}(z_0)}(u)\geq C_\delta^\prime\epsilon r,\;\;\text{for}\;\epsilon r\leq d(z_0,\partial(\Omega\setminus D)),
\label{jx0y0}
\end{equation}
for some $C_\delta^\prime>0$. We choose $r$ by imposing
\[C_\delta^\prime\epsilon r=2\tilde{C}\epsilon^\frac{1}{2}\;\;\Leftrightarrow\;\;r=\frac{2\tilde{C}}{C_\delta^\prime}\epsilon^{-\frac{1}{2}},\]
where $\tilde{C}$ is the constant in Proposition \ref{LB*}.
With this choice of $r$ we see that the existence of $z_0\in\Omega\setminus D$, $d(z_0,\partial(\Omega\setminus D))\geq \frac{2\tilde{C}}{C_\delta^\prime}\epsilon^\frac{1}{2}$, implies by Proposition \ref{LB*} and \eqref{jx0y0}
\[J_\Omega^\epsilon(u)\geq 2l\sigma+\tilde{C}\epsilon^\frac{1}{2}\]
that contradicts Proposition \ref{UB2} for small $\epsilon$. It follows
\begin{equation}
z\in\Omega\setminus D,\;d(z,\partial(\Omega\setminus D))\geq C_\delta^*\epsilon^\frac{1}{2}
\;\;\Rightarrow\;\vert u(z)-a_z\vert<\delta,\;\text{for some}\; a_z\in A,
\label{extra*}
\end{equation}
where we have set $C_\delta^*= \frac{2\tilde{C}}{C_\delta^\prime}$.

2. $a_z\equiv a_-$. The smoothness of $\partial\Omega$ implies the existence of $d>0$ such that the points $z\in C_d$, $C_d=\{z\in\Omega: d(z,\partial\Omega)\leq d\}$ can be represented with coordinates $(s,r)$ where $r=d(z,\partial\Omega)$ and $s$ is the arclength of the orthogonal projection of $z$ on $\partial\Omega$.

 Given $z_0=z(0,s_0)\in\partial\Omega$ and a small interval $(s_1,s_2)$ that contains $s_0$ set
 \[N_{s_0,(s_1,s_2)}=\{z=z(r,s): r\in(0,C_\delta^*\epsilon^\frac{1}{2}),\;s\in(s_1,s_2)\}.\]
 Assume that $d(N_{s_0,(s_1,s_2)},D)\geq C_\delta^*\epsilon^\frac{1}{2}$ and observe that $a_{z( C_\delta^*\epsilon^\frac{1}{2},s_0)}\neq a_-$ implies, via the continuity of $u$ and Lemma \ref{lower-sigma},
 \begin{equation}
 J^\epsilon_{N_{s_0,(s_1,s_2)}}\geq(\sigma-C_W\delta^2)(s_2-s_1)(1-\frac{C_\delta^*\epsilon^\frac{1}{2}}{\rho}),
 \label{Jgamma}
 \end{equation}
 where $\rho>0$ is a lower bound for the the radius of curvature of $\partial\Omega$. Since the proof of Proposition \ref{LB*} implies
 \begin{equation}
 J_{D}^\epsilon(u)\geq 2l\sigma-\tilde{C}\epsilon^\frac{1}{2},
 \label{JJD}
 \end{equation}
 \eqref{Jgamma} is in contradiction with the upper bound given by Proposition \ref{UB2} and we conclude
 \[a_{z( C_\delta^*\epsilon^\frac{1}{2},s_0)}=a_-.\]
 This, the fact that $s_0$ is arbitrary and the continuity of $u$ implies the claim and we have

\begin{equation}
z\in\Omega\setminus D,\;d(z,\partial(\Omega\setminus D))\geq C_\delta^*\epsilon^\frac{1}{2}\;\;\Rightarrow\;
\vert u(x,y)-a_-\vert\leq\delta.
\label{noupB}
\end{equation}

3. We now show that this can be upgraded to \eqref{uptoB} thus eliminating $\partial\Omega^-$, that is the part of $\partial(\Omega\setminus D)$ near which we do not expect a boundary layer, and thus allowing $z$ all the way to the boundary $\partial\Omega^-$. The tool for accomplishing this is the cut-off Lemma (Theorem 4.1 in \cite{afs}). We now describe a construction which is needed for closing the potential boundary layer. Let $C>0$ a number to be chosen later and  assume that for some $z_0=z(r_0,s_0)$ with $r\in(0,C_\delta^*\epsilon^\frac{1}{2})$ and $d(z_0,D)\geq C\epsilon^\frac{1}{2}$ it results $\vert u(z_0)-a_-\vert>\delta$. Since $\delta>0$ is small Lemma \ref{Wlemma} implies $W(u)\geq\frac{c_W^2}{2}\min\{\vert u(z(r,s_0))-a_-\vert^2,\delta^2\}$, $r\in[0,r_0]$. It follows
\begin{equation}
J_{s_0}^\epsilon(u)=:\int_0^{r_0}\Big(\frac{\epsilon}{2}\vert u^\prime\vert^2+\frac{1}{\epsilon} W(u)\Big)dr\geq
\int_0^{r_0}\frac{1}{2}(\epsilon\vert v^\prime\vert^2+\frac{c_W^2}{\epsilon}\vert v\vert^2)dr,
\label{Jbar}
\end{equation}
where we have set $v=\min\{\vert u-a_-\vert,\delta\}$.
The map $\hat{v}(r)=\delta\frac{\sinh{\frac{c_W}{\epsilon}r}}{\sinh{\frac{c_W}{\epsilon}r_0}}n$, $n=\frac{u(z_0)-a_-}{\vert u(z_0)-a_-\vert}$, minimizes $\int_0^{r_0}\frac{1}{2}(\epsilon\vert v^\prime\vert^2+\frac{c_W^2}{\epsilon}\vert v\vert^2)dr $ in the class of maps that satisfy $v(0)=0$, $v(r_0)=\delta n$ and it results $\int_0^{r_0}\frac{1}{2}(\epsilon\vert\hat{v}^\prime\vert^2+\frac{c_W^2}{\epsilon}\vert\hat{v}\vert^2)dr
=\frac{\delta^2c_W}{2}\frac{\cosh{\frac{c_W}{\epsilon}r_0}}{\sinh{\frac{c_W}{\epsilon}r_0}}$ that together with
\eqref{Jbar} imply
\begin{equation}
J_{s_0}^\epsilon(u)\geq\frac{c_W\delta^2}{2}.
\label{JbarLB}
\end{equation}
Let $S=\{s_0: \vert u(z(r_0,s_0))-a_-\vert\geq\delta,\;r_0\in(0,C_\delta^*\epsilon^\frac{1}{2}),\;
d(z(r_0,s_0),D)\geq C\epsilon^\frac{1}{2}\}$ and choose $C>0$ to ensure that the segment with extreme $z(0,s_0)$ and $z(C_\delta^*\epsilon^\frac{1}{2},s_0)$ is contained in $\Omega\setminus D$. Then \eqref{JbarLB} implies
\begin{equation}
J^\epsilon_{\Omega\setminus D}(u)\geq\int_SJ_{s_0}^\epsilon(u)(1-\frac{C_\delta^*\epsilon^\frac{1}{2}}{\rho})d s_0
\geq\mathcal{H}^1(S)\frac{c_W\delta^2}{2}(1-\frac{C_\delta^*\epsilon^\frac{1}{2}}{\rho}),
\label{JbarS}
\end{equation}
This, \eqref{JbarS}, \eqref{JJD}  and \eqref{UB*} imply
\[\mathcal{H}^1(S)\frac{c_W\delta^2}{2}(1-\frac{C_\delta^*\epsilon^\frac{1}{2}}{\rho})+ 2l\sigma-\tilde{C}\epsilon^\frac{1}{2}\leq
2l\sigma +C_1\epsilon\vert\ln{\epsilon}\vert^3\] and we obtain
\begin{equation}
\mathcal{H}^1(S)\leq2\frac{\tilde{C}\epsilon^\frac{1}{2}+C_1\vert\epsilon\ln{\epsilon}\vert^3}{c_W\delta^2
(1-\frac{C_\delta^*\epsilon^\frac{1}{2}}{\rho})}
\leq\hat{C}\epsilon^\frac{1}{2}.
\label{barS-bound}
\end{equation}
 From \eqref{barS-bound} it follows that, given $z(0,s_0)\in\partial\Omega^-$ with $d(z(0,s_0),D)\geq C\epsilon^\frac{1}{2}$ for some $C>0$ sufficiently large, there are $s_1<s_0<s_2$ with $s_2-s_1\leq 2\hat{C}\epsilon^\frac{1}{2}$ such that (recalling also that we have $u=a_-$ on $\partial\Omega^-$)
 \[\vert u(z)-a_-\vert\leq\delta,\;\;z\in\partial N_{s_0,(s_1,s_2)}.\]
 This, since $\delta>0$ is a small number, implies that we can invoke Theorem 4.1 in \cite{afs} and conclude that \[\vert u(z)-a_-\vert\leq\delta,\;\;z\in N_{s_0,(s_1,s_2)}.\]
 This and the fact that the only condition imposed to $s_0$ is $d(z(0,s_0),D)\geq C\epsilon^\frac{1}{2}$ imply \eqref{uptoB} for some constant $C_\delta>0$.

4. Finally we establish \eqref{exp-refin} by applying linear theory. Set $\Omega^\epsilon=\{\zeta: \epsilon\zeta\in\Omega\}$.
 The map $U^\epsilon:\Omega^\epsilon\rightarrow\R^m$, $U^\epsilon(\zeta)=u(\epsilon\zeta)$ is a minimizer of $\int_{\Omega^\epsilon}\Big(\frac{1}{2}\vert\nabla V\vert^2+W(V)\Big)d\xi d\eta$. Fix $z\in\Omega$ with $d(z,D)\geq C\epsilon^\frac{1}{2}$. Now \eqref{uptoB} implies that we can assume $C>0$ sufficiently large to ensure
 \[\vert u-a_-\vert\leq\delta,\quad\text{on}\quad B_R(z)\cap\Omega,\quad R=d(z,D)-C\epsilon^\frac{1}{2}.\]
 Set $E=B_R(z)\cap\Omega$ and $E^\epsilon=\{\zeta:\epsilon\zeta\in E\}$. The minimizer $U^\epsilon$ satisfies
 \begin{equation}
 \begin{split}
 &\vert U^\epsilon-a_-\vert\leq\delta,\quad\text{on}\quad E^\epsilon,\\
 & U^\epsilon=a_-,\quad\text{on}\quad \partial E^\epsilon\cap\partial\Omega^\epsilon.
 \end{split}
 \label{B-cond-U}
 \end{equation}
 Arguing as in the proof of  Lemma 4.4 (pag.123 in \cite{afs}) with $A=E^\epsilon$ and $r=\delta$ we deduce
 \[\vert U^\epsilon(\frac{z}{\epsilon})-a_-\vert^2\leq\delta^2\varphi(\frac{z}{\epsilon}),
\]
where $\varphi:E^\epsilon\rightarrow\R$ is the solution of
\begin{equation}
\begin{split}
&\Delta\varphi=c_W^2\varphi,\quad\text{on}\quad E^\epsilon,\\
&\varphi=\frac{1}{\delta^2}\vert U^\epsilon-a_-\vert^2\leq 1,\quad\text{on}\quad \partial E^\epsilon.
\end{split}
\label{varphi}
\end{equation}
Since from \eqref{B-cond-U} we have $\varphi=0$ on $\partial E^\epsilon\cap\partial\Omega^\epsilon$, the restriction to $ E^\epsilon$ of the solution $\psi$ of the problem
\begin{equation}
\begin{split}
&\Delta\psi=c_W^2\psi,\quad\text{on}\quad B_{\frac{R}{\epsilon}}(\frac{z}{\epsilon}),\\
&\varphi=1,\quad\text{on}\quad\partial B_{\frac{R}{\epsilon}}(\frac{z}{\epsilon}),
\end{split}
\label{varpsi}
\end{equation}
is a supersolution for problem \eqref{varphi} therefore, using also Lemma A.1 in  \cite{afs}, we have
 \begin{equation}
 \vert U^\epsilon(\frac{z}{\epsilon})-a_-\vert^2\leq\delta^2\psi(\frac{z}{\epsilon})
 \leq\delta^2e^{-k_0\frac{R}{\epsilon}},\;\;\;\text{for}\;\;\;R\geq\epsilon\rho_0,
 \label{U-exp-est}
 \end{equation}
for appropriate $k_0>0$ and $\rho_0$. It follows that
\[\vert u(z)-a_-\vert\leq\delta e^{-\frac{k_0}{2\epsilon}(d(z,D)-C\epsilon^\frac{1}{2})},\;\;\;\text{for}\;\;\;d(z,D)\geq C\epsilon^\frac{1}{2},\]
where, by adjusting the constant $C>0$ we have absorbed the term $\epsilon\rho_0$ in $C\epsilon^\frac{1}{2}$.
 The proof is complete.
\end{proof}
A first consequence of Lemma \ref{nearone*} is a refinement of the lower bound \eqref{LB*1}. Set $\hat{y}_1=\frac{1}{2}(y_1+y_2)-\frac{1}{4}(y_2-y_1)$, $\hat{y}_2=\frac{1}{2}(y_1+y_2)+\frac{1}{4}(y_2-y_1)$, then, from $y_2-y_1\geq\frac{1}{4}(h-l)$ and the exponential estimate in Lemma \ref{nearone*} it follows
\[\begin{split}
&x\in(C_0\epsilon,l-C_0\epsilon),\;j=1,2\quad\Rightarrow\\
&\vert u(x,\hat{y}_j)-a_-\vert\leq\delta e^{-\frac{k}{16\epsilon}(h-l-C_\delta\epsilon^\frac{1}{2})}\leq\delta e^{-\frac{k}{32\epsilon}(h-l)}.
\end{split}\]
Then,  by taking $\epsilon>0$ sufficiently small we can assume

\[\vert u(x,\hat{y}_j)-a_-\vert\leq\epsilon^\frac{1}{2},\;\;x\in(C_0\epsilon,l-C_0\epsilon),\;j=1,2.\]
Therefore proceeding  as for \eqref{halph} with $\hat{y}_j$ in place of $y_j$ and $\delta=\epsilon^\frac{1}{2}$ we obtain
\begin{equation}
 J^\epsilon_{(0,\hat{y}_1)}(u(x,\cdot)), J^\epsilon_{(\hat{y}_2,h)}(u(x,\cdot))
\geq\sigma-C_W\epsilon,\;\;x\in(C_0\epsilon,l-C_0\epsilon),
\label{small-x}
\end{equation}
and therefore the lower bound
\begin{equation}
J_\Omega^\epsilon(u)\geq\int_0^lJ^\epsilon_{(0,h)}(u(x,\cdot))dx\geq 2l\sigma-\tilde{C}^\prime\epsilon,\;\;\;\tilde{C}^\prime=4\sigma C_0+C_Wl.
\label{LB-ref}
\end{equation}

The estimate \eqref{LB-ref} is basic for deriving information on the structure of the minimizer $u$ in the set $[0,l]\times([0,y_1]\cup[y_2,h])$. We begin by studying the one-dimensional problem:
\subsubsection{The one-dimensional problem.}\label{ONE}
Consider the minimization problem
\begin{equation}
\begin{split}
&\min_{v\in\mathscr{V}}J^\epsilon_{(0,\lambda)}(v),\;\;J^\epsilon_{(0,\lambda)}(v)=\int_0^\lambda(\frac{\epsilon}{2}\vert v^\prime\vert^2+\frac{1}{\epsilon}W(v))ds,\\
&\mathscr{V}=\{v\in H^1((0,\lambda);\R^m):v(0)=a_-+v_0,\;\vert v_0\vert<\epsilon^\frac{1}{2},\;v(\lambda)=a_+\}
\end{split}
\label{problem*2}
\end{equation}
\begin{lemma}
 \label{min-0,lambda}
 Assume $v$ is a minimizer of \eqref{problem*2}. Then there exist $0<s^-<s^+<\lambda$ and a constant $C^*>0$ such that
\begin{equation}
\begin{split}
&s^+-s^-\leq C^*\epsilon^\frac{1}{2},\;\;(C^*=\frac{4\sigma}{c_W^2})\\
&\text{and}\\
&v(s)\in B_{\epsilon^\frac{1}{4}}(a_-),\;\;s\in[0,s^-),\\
&v(s)\in\R^m\setminus\cup_{a\in A}B_{\epsilon^\frac{1}{4}}(a),\;\;s\in[s^-,s^+],\\
&v(s)\in B_{\epsilon^\frac{1}{4}}(a_+),\;\;s\in(s^+,\lambda).
\end{split}
\label{properties}
\end{equation}

\end{lemma}
\begin{proof}
The same argument that establishes \eqref{jbarv} implies that there is a constant $\bar{C}>0$ such that
\begin{equation}
J^\epsilon_{(0,\lambda)}(v)\leq\sigma+\bar{C}\epsilon.
\label{upper-v}
\end{equation}
Set
\[\begin{split}
& s^-=\sup\{s>0: v(s)\in B_{\epsilon^\frac{1}{4}}(a_-)\},\\
& s^+=\inf\{s<\lambda: v(s)\in B_{\epsilon^\frac{1}{4}}(a_+)\}.
\end{split}\]
Lemma 2.4 and Lemma 2.5 in \cite{afs} imply that, for $\delta>0$ small and $a\in A$, the existence of $s_1<s_2<s_3$ such that
\[v(s_i)\in B_\delta(a),\;i=1,3\quad\text{and}\;v(s_2)\not\in B_\delta(a)\]
is incompatible with the minimality of $v$.
This property of minimizers implies that \eqref{properties}$_2$ and
\eqref{properties}$_4$ hold and moreover that
\begin{equation}
v(s)\in\R^m\setminus (B_{\epsilon^\frac{1}{4}}(a_-)\cup B_{\epsilon^\frac{1}{4}}(a_+)),\;\;s\in[s^-,s^+].
\label{not-apm}
\end{equation}
From the minimality of $\bar{u}$ and the upper bound \eqref{upper-v} we also have
\[v(s)\not\in B_{\epsilon^\frac{1}{4}}(a),\;\;a\in A\setminus\{a_-,a_+\},\;s\in[s^-,s^+].\]
From this and \eqref{not-apm}
 \eqref{properties}$_3$ follows. From \eqref{properties}$_3$ and Lemma \ref{Wlemma} we have $W(v(s))\geq \frac{1}{2}c_W^2\epsilon^\frac{1}{2}$ and therefore from \eqref{upper-v} we obtain
\[(s^+-s^-)\frac{1}{2\epsilon}c_W^2\epsilon^\frac{1}{2}\leq\sigma+\bar{C}\epsilon\]
and \eqref{properties}$_1$ follows with $C^*=\frac{4\sigma}{c_W^2}$. The proof is complete.
\end{proof}
Define $\mathscr{V}_a=\mathscr{V}_a^1\cup\mathscr{V}_a^2$ where
\begin{equation}
\begin{split}
&\mathscr{V}_a^1=\{v\in\mathscr{V}:\,v(\bar{s})\in B_{\epsilon^\frac{1}{4}}(a),\;\;\text{for some}\;\;a\neq a_\pm,\,\text{and}\;\bar{s}\in(0,\lambda)\},\\
&\mathscr{V}_a^2=\{v\in\mathscr{V}:\,v(\bar{s}_i)\in B_{\epsilon^\frac{1}{4}}(a_i),\;a_1=a_-,\,a_2=a_+,\,a_3=a_-\\
&\text{for some}\;0\leq\bar{s}_1<\bar{s}_2<\bar{s}_3<\lambda\}.
\end{split}
\label{Wa}
\end{equation}
Let $\mathscr{W}=\mathscr{V}\setminus\mathscr{V}_a$.
For $w\in\mathscr{W}$ we set
 \[\begin{split}
 &s^{-,w}=\sup\{s\in(0,\lambda):w(s)\in B_{\epsilon^\frac{1}{4}}(a_-),\\
 &s^{+,w}=\inf\{s\in(0,\lambda):w(s)\in B_{\epsilon^\frac{1}{4}}(a_+),\\
 &S^w=\{s\in(0,\lambda): w(s)\in\R^m\setminus\cup_{a\in A}B_{\epsilon^\frac{1}{4}}(a)\}.
 \end{split}\]
 Let $K>1$ a number to be chosen later and define
 \begin{equation}
 \begin{split}
 &\mathscr{W}^*=\Big\{
w\in\mathscr{W}:\,\vert S^w\vert < 2C^*\epsilon^\frac{1}{2}
,\;\text{and}\\ &\left.\begin{array}{l}
w(s)\in B_{K\epsilon^\frac{1}{4}}(a_-),\;s\in[0,s^{-,w}),\\
w(s)\in B_{K\epsilon^\frac{1}{4}}(a_+),\;s\in(s^{+,w},\lambda]
\end{array}.
\right.\Big\}\end{split}
 \label{v-fat}
 \end{equation}
 Note that
  \begin{equation}
 s^{+,w}-s^{-,w}\leq 2C^*\epsilon^\frac{1}{2},\;\;w\in\mathscr{W}^*,
 \label{sw-sw}
 \end{equation}
 and that Lemma \ref{lower-sigma} implies

 \begin{equation}
 J^\epsilon_{(0,\lambda)}(w)\geq\sigma-C_W\epsilon^\frac{1}{2},\;\;w\in\mathscr{W}^*.
 \label{JW*}
 \end{equation}
 Note also that the maps that satisfies the \eqref{properties}, in particular minimizers of \eqref{problem*2}, belong to $\mathscr{W}^*$.

 We divide $\mathscr{W}^c=\mathscr{W}\setminus\mathscr{W}^*$  in  two disjoint parts $\hat{\mathscr{W}}^c$ and $\tilde{\mathscr{W}}^c$ defined as follows

\begin{equation}
\begin{split}
&\hat{\mathscr{W}}^c=\{w\in\mathscr{W}^c:
\,\vert S^w\vert\geq 2C^*\epsilon^\frac{1}{2}\},\\
&\tilde{\mathscr{W}}^c=\{w\in\mathscr{W}^c:
\,\vert S^w\vert< 2C^*\epsilon^\frac{1}{2},\;\text{and}\\
&\;w(\bar{s})\not\in B_{K\epsilon^\frac{1}{4}}(a_-)\cup B_{K\epsilon^\frac{1}{4}}(a_+),
\;\text{for some}\;\bar{s}\not\in[s^{-,w},s^{+,w}]\}.
\end{split}
\label{fourW}
\end{equation}

We set $\mathscr{V}^c=\mathscr{V}_a\cup\hat{\mathscr{W}}^c\cup\tilde{\mathscr{W}}^c$ and observe that $\mathscr{V}^c=\mathscr{V}\setminus\mathscr{W}^*$. From the definition of $\mathscr{V}^c$ we see that maps in
$\mathscr{V}^c$ have a structure that differs substantially from that of minimizers of problem \eqref{problem*2}
 described in Lemma \ref{min-0,lambda}. We need a quantitative estimate of the energy price paid by a map $w\in\mathscr{V}^c$.
 \begin{lemma}
 \label{winVc}
 Let $\mathscr{V}^c$ be as before. Then there is $K>1$ such that
 \begin{equation}
 w\in\mathscr{V}^c\;\;\Rightarrow\;\;J^\epsilon_{(0,\lambda)}(w)\geq\sigma+C_W\epsilon^\frac{1}{2}.
\label{oned-est}
 \end{equation}
 \end{lemma}
 \begin{proof}

 1. If $w\in\mathscr{V}_a^1$ we have $J^\epsilon_{(0,\lambda)}(w)\geq\sigma+\eta$ for some $\eta>0$. This follows from the minimality of $\bar{u}$. If $w\in\mathscr{V}_a^2$ Lemma \ref{lower-sigma} yields $J^\epsilon_{(0,\lambda)}(w)\geq 2\sigma-C_W\epsilon^\frac{1}{2}$.

 2. If $w\in\hat{\mathscr{W}}^c$ we have
 \[W(w(s))\geq\frac{1}{2}c_W^2\epsilon^\frac{1}{2},\;\;s\in S^w.\]
 it follows
 \[J^\epsilon_{(0,\lambda)}(w)\geq\frac{1}{\epsilon}W(w(s))\vert S^w\vert\geq 2c_W^2C^*=4\sigma,\;\;(C^*=\frac{2\sigma}{c_W^2}).\]

 3. If $w\in\tilde{\mathscr{W}}^c$ we have
 \[\begin{split}
 &w(s^{-,w})\in\bar{B}_{\epsilon^\frac{1}{4}}(a_-),\\
 &w(s^{+,w})\in\bar{B}_{\epsilon^\frac{1}{4}}(a_+).
 \end{split}\]
 This and Lemma \ref{lower-sigma} imply
 \begin{equation}
 J^\epsilon_{(s^{-,w},s^{+,w})}(w)\geq\sigma-C_W\epsilon^\frac{1}{2}.
 \label{partialE}
 \end{equation}
By assumption there is $\bar{s}\not\in(s^{-,w},s^{+,w})$ that satisfies
 $w(\bar{s})\in\R^m\setminus(B_{K\epsilon^\frac{1}{4}}(a_-)\cup B_{K\epsilon^\frac{1}{4}}(a_+))$. From this it  follows the existence of two intervals $I_j=[r_j,s_j]\subset(0,\lambda)\setminus(s^{-,w},s^{+,w})$, $j=1,2$ such that
 \[\begin{split}
 &W(w(s))\geq\frac{1}{2}c_W^2\epsilon^\frac{1}{2},\;\,s\in I_j,\\
 &\vert\vert w(r_j)\vert-\vert w(s_j)\vert\vert\geq (K-1)\epsilon^\frac{1}{4}.
 \end{split}\]
 This and a standard computation imply
 \[J^\epsilon_{I_j}(w)\geq\frac{1}{2}\Big(\epsilon^\frac{3}{2}\frac{(K-1)^2}{\vert I_j\vert}
 +\frac{1}{\epsilon^\frac{1}{2}}c_W^2\vert I_j\vert\Big)\geq c_W(K-1)\epsilon^\frac{1}{2},\]
 that together with \eqref{partialE} yields
\[J^\epsilon_{(0,\lambda)}(w)\geq\sigma+(2c_W(K-1)-C_W)\epsilon^\frac{1}{2}\geq\sigma+C_W\epsilon^\frac{1}{2},\]
where we have chosen $K=1+\frac{C_W}{c_W}$. The proof is complete.
 \end{proof}

\subsubsection{The structure inside the domain}
\label{Structure}
From \eqref{small-x} we have a good lower bound for the energy $J_{(0,\hat{y}_1)}(u(x,\cdot))$ or $J_{(\hat{y}_2,h)}(u(x,\cdot))$ of the restrictions of  the minimizer $u:\Omega\rightarrow\R^m$ to each fiber $\{x\}\times(0,\hat{y}_1)$, or $\{x\}\times(\hat{y}_2,h)$ and therefore we expect that  $u(x,\cdot)\vert_{(0,\hat{y}_1)}$ and $u(x,\cdot)\vert_{(\hat{y}_2,h)}$ should be well approximated by a translation of $\bar{u}(-\cdot)$ and $\bar{u}$ respectively. We prove that this is indeed the case for most $x\in(C_0\epsilon,l-C_0\epsilon)$. To show this we apply the results in Lemma \ref{min-0,lambda} and Lemma \ref{winVc} to $u$ restricted to the fibers $\{x\}\times(0,\hat{y}_1)$ and $\{x\}\times(\hat{y}_2,h)$. This requires a natural reinterpretation of the notation. For instance, when dealing with the fiber $\{x\}\times(\hat{y}_2,h)$ the interval $(0,\lambda)$ in problem \eqref{problem*2} should be replaced by the interval $(\hat{y}_2,h)$ and the intervals $[0,s^{-,w}]$ and $[s^{+,w},\lambda]$ in the definition \eqref{v-fat} of $\mathscr{W}^*$ with the intervals $[\hat{y}_2,\hat{y}_2+s^{-,w}]$ and $[\hat{y}_2+s^{+,w},h]$. For $w=u(x,\cdot)\vert_{(\hat{y}_2,h)}$ we set
\[
\eta(x)=h-\hat{y}_2-s^{+,u(x,\cdot)\vert_{(\hat{y}_2,h)}}.
\]
Let $X=X_1\cup X_2\subset(C_0\epsilon,l-C_0\epsilon)$ be defined by
\[\begin{split}
& X_1=\{x\in(C_0\epsilon,l-C_0\epsilon): u(x,\cdot)\vert_{(0,\hat{y}_1)}\in\mathscr{V}^c\},\\
& X_2=\{x\in(C_0\epsilon,l-C_0\epsilon):
u(x,\cdot)\vert_{(\hat{y}_2,h)}\in\mathscr{V}^c\}.
\end{split}\]
\begin{lemma}
\label{bad-set}
There exists a constant $C^\prime>0$ such that
\begin{equation}
\mathcal{H}^1(X)+\leq C^\prime\epsilon^\frac{1}{2}\vert\ln{\epsilon}\vert^3.
\label{X-bound}
\end{equation}
\end{lemma}
\begin{proof}
From Lemma \ref{winVc}, \eqref{small-x} and Proposition \ref{UB2} we have
\[\sum_{i\in\{1,2\}}\Big((\sigma-C_W\epsilon)(l-C_0\epsilon-\mathcal{H}^1(X_i))
+\mathcal{H}^1(X_i)(\sigma+C_W\epsilon^\frac{1}{2})\Big)\leq J_\Omega^\epsilon(u)\leq2l\sigma +
C_1\epsilon\vert\ln{\epsilon}\vert^3.\]
It follows $\mathcal{H}^1(X_1)+\mathcal{H}^1(X_2)\leq C^\prime\epsilon^\frac{1}{2}\vert\ln{\epsilon}\vert^3$
 with $C^\prime=\frac{C_1}{C_W}$. This and $\mathcal{H}^1(X_1\cup X_2)\leq\mathcal{H}^1(X_1)+\mathcal{H}^1(X_2)$
 conclude the proof.
\end{proof}
\vskip.2cm
We now continue with the proof of Theorem \ref{BL-exist}.
\begin{proof} (of Theorem \ref{BL-exist})

1. For $\xi\in(0,\frac{l}{2}]$ we let $\Omega(\xi)$ be the connected component of $\Omega\cap(-2C_\delta\epsilon^\frac{1}{2},\xi)\times(\hat{y}_2,+\infty)$ that contains $(0,\xi)\times(\hat{y}_2,h)$ (see Figure \ref{alpha}).

\begin{figure}
  \begin{center}
\begin{tikzpicture}[scale=1.2]
\draw [blue] (-3.3,0) -- (-1,0);

\draw [] (-3,2) -- (-1,2);
\draw [] (1,2) -- (3,2);

\draw [] (-1,2) -- (-1,0);
\draw [] (3,0) -- (3,2);

\draw[blue] (-3.3,1.95) to [out=15,in=180] (-3,2);
\draw[] (.7,2.05) to [out=-15,in=180] (1,2);

\draw [blue] (-3.3,1.95) -- (-3.3,0);
\draw [] (.7,2.05) -- (.7,0);

\path[fill=lightgray] (2.4595,0)--(2.5946,0)--(3,1.5)--(3,2)--(2.4595,0);
\draw []  (2.4595,0)--(2.5946,0)--(3,1.5)--(3,2)--(2.4595,0);

\path[fill=lightgray] (.7,.3348)--(.7,0)--(2.4595,0)--(3,2)--(1.15,2)--(.7,.3348);
\draw [] (.7,.3348)--(.7,0)--(2.4595,0)--(3,2)--(1.15,2)--(.7,.3348);

\draw [red] (-2.85,2) -- (-1,2);
\draw [red] (1.15,2) -- (3,2);

\draw [red] (-1,1.5) -- (-1,2);
\draw [red] (3,1.5) -- (3,2);

\draw [dotted] (-1,1.5) -- (-2.85,2);
\draw [dotted] (3,1.5) -- (1.15,2);

\draw [dotted] (-3.,2) -- (-3.,0);
\draw [dotted] (1,2.) -- (1,0);

\draw [blue] (.7,0) -- (3,0);
\draw [blue](.7,2) -- (.7,0);

\node [below] at (1,-.5) {$(0,\hat{y}_2)$};
\draw [->](1,-.5) -- (1,-.1);

\node [above] at (0,2.5) {$\Omega(\xi)$};
\draw [->] (0,2.5) -- (1.5,2.1);
\draw [->](0,2.5) -- (-1.5,2.1);

\node [right] at (2.85,2.2) {$(\xi,h)$};
\node [right] at (2.90,1.5) {$(\xi,h-\eta(\xi))$};
\node [right] at (2.85,-.2) {$(\xi,\hat{y}_2)$};

\draw [blue](.65,2)--(.75,2);
\node [left] at (.7,1.5) {$H$};
\node [below] at (1.9,0) {$L_\xi$};

\draw [->](-3.7,2.35)--(-3.35,1.5);
\draw [->](-3.7,2.35)--(-3.16,1.99);
\node [left] at (-3.7,2.35) {$H$};

\draw [->](-1,-.5)--(.65,-.1);
\node [below] at (-1,-.5) {$(-2C_\delta\epsilon^\frac{1}{2},\hat{y}_2)$};

\node [] at (1.8,.8) {$\Omega(\xi)^\prime$};

\draw [->](3.7,.7)--(2.7,.7);
\node [right] at (3.7,.7) {$\Omega(\xi)^{\prime\prime}$};

\draw [] (-2.2,2) arc [radius=.7, start angle=0, end angle= -15.123];

\node [above] at (-2.3,2.1) {$\alpha_\xi$};
\draw[->] (-2.3,2.2)--(-2.3,1.95);

\end{tikzpicture}
\end{center}
\caption{$\Omega(\xi), \Omega(\xi)^\prime, \Omega(\xi)^{\prime\prime}, L_\xi$ and $H$ for two different $\Omega$.}
\label{alpha}
\end{figure}
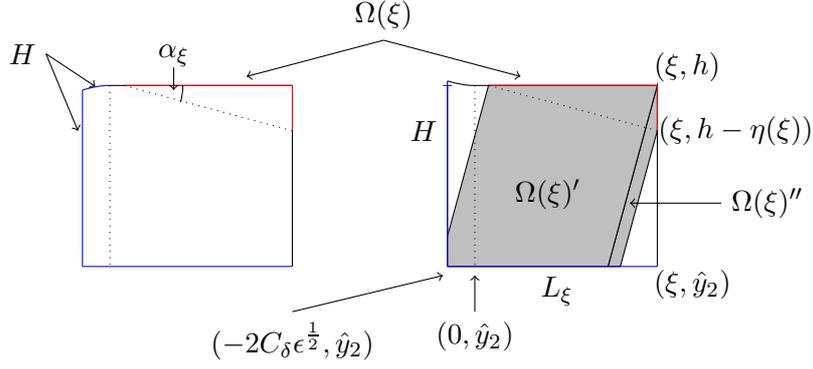

Set
\[\begin{split}
&L_\xi=[-2C_\delta\epsilon^\frac{1}{2},\xi]\times\{\hat{y}_2\},\\
&H=\partial\Omega(\xi)\cap(-\infty,\xi)\times(\hat{y}_2,h),
\end{split}\]
and observe that, for small $\epsilon>0$, the estimate \eqref{exp-refin} in Lemma \ref{nearone*} implies
\begin{equation}
\vert u(z)-a_-\vert<\epsilon,\;\;\text{on}\;\;L_\xi\cup H,
\label{B-cond}
\end{equation}

Also note that, if $\xi\in[C_0\epsilon,\frac{l}{2}]\setminus X_2$ we have $u(x,\cdot)\vert_{(\hat{y}_2,h)}\in\mathscr{W}^*$ and
therefore
\begin{equation}
\vert u(z)-a_+\vert<K\epsilon^\frac{1}{4},\;\;\text{on}\;\;H_\eta,
\label{B-cond1}
\end{equation}
where
\[H_\eta=\{\xi\}\times[h-\eta(\xi),h].\]

\vskip.2cm
2. Assume that $\xi\in[C_0\epsilon,\frac{l}{2}]\setminus X_2$ and set

\[\sin{\alpha_{\xi}}=\frac{\eta(\xi)}{(\xi^2+\eta^2(\xi))^\frac{1}{2}},\quad \cos{\alpha_{\xi}}=\frac{\xi}{(\xi^2+\eta^2(\xi))^\frac{1}{2}}.\]

We regard $\Omega(\xi)$ as a union of fibers parallel to $(\sin{\alpha_{\xi}},\cos{\alpha_{\xi}})$ and let $\Omega(\xi)^\prime$ be the union of the fibers that terminate on $[C_0\epsilon,\xi]\times\{h\}$ and
$\Omega(\xi)^{\prime\prime}$ the union of the fibers that terminate on $H_\eta$. Since $u=a_+$ on $[C_0\epsilon,\xi]\times\{h\}$ and each fiber with the second extreme on $[C_0\epsilon,\xi]\times\{h\}$ has its first extreme on $L_\xi\cup H$. from \eqref{B-cond} it follows

\begin{equation}
J^\epsilon_{\Omega(\xi)^\prime}(u)\geq(\sigma-C_W\epsilon)(\xi-C_0\epsilon-\mathcal{H}^1(X))\cos{\alpha_{\xi}}
\geq\sigma\xi\cos{\alpha_{\xi}}-C\epsilon^\frac{1}{2}\vert\ln{\epsilon}\vert^3,
\label{Js1}
\end{equation}
where we have used \eqref{X-bound} and $C>0$ denotes a generic constant independent of $\epsilon$.

On the other hand  \eqref{B-cond1} implies that the  contribution $J_{\Omega(\xi)^{\prime\prime}}(u)$ of the fibers that intersect $H_\eta$ satisfies

\[J^\epsilon_{\Omega(\xi)^{\prime\prime}}(u)\geq(\sigma-K^2C_W\epsilon^\frac{1}{2}))\eta(\xi)\sin{\alpha_{\xi}}
\geq\sigma\eta(\xi)\sin{\alpha_\xi}-C\epsilon^\frac{1}{2}.\]

This and \eqref{Js1} yield
\begin{equation}
J^\epsilon_{\Omega(\xi)}\geq J^\epsilon_{\Omega(\xi)^\prime}(u)+J^\epsilon_{\Omega(\xi)^{\prime\prime}}(u)\geq
\sigma\eta(\xi)\sin{\alpha_\xi}+\sigma\xi\cos{\alpha_{\xi}}-C\epsilon^\frac{1}{2}\vert\ln{\epsilon}\vert^3.
\label{Js}
\end{equation}

From \eqref{small-x} with $R=(C_0\epsilon,l-C_0\epsilon)\times(0,\hat{y}_1)$ and $R(\xi)=(\xi,l-C_0\epsilon)\times(\hat{y}_2,h)$ we also have
\begin{equation}
\begin{split}
&J_R^\epsilon(u)\geq(\sigma-C_W\epsilon)(l-2C_0\epsilon)\geq\sigma l-C\epsilon,\\
&J^\epsilon_{R(\xi)}(u)\geq(\sigma-C_W\epsilon)(l-\xi-C_0\epsilon)\geq\sigma(l-\xi)-C\epsilon.
\end{split}
\label{JR}
\end{equation}

The upper bound \eqref{UB*} and the estimate \eqref{Js} ,\eqref{JR} yield

\begin{equation}
\eta(\xi)\sin{\alpha_\xi}\leq\xi(1-\cos{\alpha_\xi})+C\epsilon^\frac{1}{2}\vert\ln{\epsilon}\vert^3.
\label{eta-cond}
\end{equation}
The definition of $\sin{\alpha_\xi}$ and $\cos{\alpha_\xi}$ implies that \eqref{eta-cond} is equivalent to
\[\eta(\xi)^2\leq 2C\xi\epsilon^\frac{1}{2}\vert\ln{\epsilon}\vert^3+C^2\epsilon\ln{\epsilon}\vert^6.\]
Since we have $\xi\leq\frac{l}{2}$ and $\xi\not\in X$ we conclude
\begin{equation}
\eta(\xi)\leq C^\sharp\epsilon^\frac{1}{4}\vert\ln{\epsilon}\vert^\frac{3}{2},\;\;\xi\in(C_0\epsilon,\frac{l}{2}]\setminus X.
\label{Eta-Bound}
\end{equation}
for some constant $C^\sharp>0$ independent of $\epsilon$.

\vskip.2cm
3.
We can perform a similar analysis for estimating $\eta(\xi)$ for $\xi\in[\frac{l}{2},l-C_0\epsilon)$ and for estimating the size of $s^{-,u(x,\cdot)\vert(0,\hat{y}_1)}$ for $\xi\in(C_0\epsilon,l-C_0\epsilon)$. Proceeding in this way and recalling the definition of $\eta(\xi)$ and the estimate \eqref{sw-sw} we finally obtain
\begin{equation}
\left.\begin{array}{l}
 h-\hat{y}_2-s^{-,u(x,\cdot)\vert_{(\hat{y}_2,h)}}\leq 2C^\sharp\epsilon^\frac{1}{4}\vert\ln{\epsilon}\vert^\frac{3}{2},\\\\
s^{+,u(x,\cdot)\vert_{(0,\hat{y}_1)}}\leq 2C^\sharp\epsilon^\frac{1}{4}\vert\ln{\epsilon}\vert^\frac{3}{2},
\end{array}\right.\;\;x\in(C_0\epsilon,l-C_0\epsilon)\setminus X.
\label{s+u-}
\end{equation}
This implies that most of the energy of the minimizer $u$ is concentrated near $\partial^+\Omega$.  Set
\[D_\epsilon=\Big((C_0\epsilon,l-C_0\epsilon)\setminus X\Big)\times\Big((0,2C^\sharp\epsilon^\frac{1}{4}\vert\ln{\epsilon}\vert^\frac{3}{2})\cup
(h-2C^\sharp\epsilon^\frac{1}{4}\vert\ln{\epsilon}\vert^\frac{3}{2},h)\Big).\]

From \eqref{s+u-} and the fact that $u(x,\cdot)\vert_{(0,\hat{y}_1)}$ and $u(x,\cdot)\vert_{(\hat{y}_2,h)}$ belong to $\mathscr{W}^*$ it follows
\begin{equation}
\vert u(z)-a_-\vert\leq K\epsilon^\frac{1}{4},\quad z\in\Big((C_0\epsilon,l-C_0\epsilon)\setminus X\Big)\times
(2C^\sharp\epsilon^\frac{1}{4}\vert\ln{\epsilon}\vert^\frac{3}{2},h-2C^\sharp\epsilon^\frac{1}{4}\vert\ln{\epsilon}\vert^\frac{3}{2})
,
\label{NEARa-}
\end{equation}
and therefore Lemma \ref{lower-sigma} and \eqref{X-bound} imply
\[J_{D_\epsilon}^\epsilon(u)\geq 2(\sigma-K^2C_W\epsilon^\frac{1}{2})(l-2C_0\epsilon-C^\prime\epsilon^\frac{1}{2}\vert\ln{\epsilon}\vert^3)
\geq 2\sigma l-C\epsilon^\frac{1}{2}\vert\ln{\epsilon}\vert^3.\]

4. Conclusion.
As in the proof of Lemma \ref{nearone*} we find that
$\min_{a\in A}\vert u(z_0)-a\vert\geq\delta$ for some $z_0\in\Omega\setminus D_\epsilon$ implies
\begin{equation}
J^\epsilon_{B_{\epsilon r}(z_0)}(u)\geq C_\delta^\prime\epsilon r,\;\;\text{for}\;\epsilon r\leq d(z_0,\partial (\Omega\setminus D_\epsilon)),
\label{jx0y0S}
\end{equation}
for some $C_\delta^\prime>0$. We choose $r$ by imposing
\[C_\delta^\prime\epsilon r=3\sigma C^\prime\epsilon^\frac{1}{2}\vert\ln{\epsilon}\vert^3\;\;\Leftrightarrow\;\;r=\frac{3 \sigma C^\prime}{C_\delta^\prime}\epsilon^{-\frac{1}{2}}\vert\ln{\epsilon}\vert^3.\]
With this choice of $r$ we see that the existence of $z_0\in\Omega\setminus D_\epsilon$, $d(z_0,\partial (\Omega\setminus D_\epsilon))\geq \frac{3\sigma C^\prime}{C_\delta^\prime}\epsilon^\frac{1}{2}\vert\ln{\epsilon}\vert^3$, implies
\[J^\epsilon_{D_\epsilon}(u)+J^\epsilon_{B_{\epsilon r}(z_0)}(u)\geq 2l\sigma+
\sigma C^\prime\epsilon^\frac{1}{2}\vert\ln{\epsilon}\vert^3.\]
This estimate, if $\epsilon>0$ is small, collides with the upper bound \eqref{UB*} and we conclude that
\begin{equation}
\min_{a\in A}\vert u(z)-a\vert<\delta,\;\;z\in\Omega\setminus D_\epsilon,\;d(z,\partial (\Omega\setminus D_\epsilon))\geq
C\epsilon^\frac{1}{2}\vert\ln{\epsilon}\vert^3,
\label{a?}
\end{equation}
where $C=\frac{3\sigma C^\prime}{C_\delta^\prime}$.
This, \eqref{NEARa-}, \eqref{uptoB} and the continuity of $u$ imply that actually we have $a=a_-$ in \eqref{a?}.
From the definition of $D_\epsilon$ it follows that for  $C^1>0$ sufficiently large
  we have, for small $\epsilon>0$,
\[\begin{split}
&z\in\Omega,\;\;d(z,\partial^+\Omega)\geq C^1\epsilon^\frac{1}{4}\vert\ln{\epsilon}\vert^\frac{3}{2}\;\;\Rightarrow\\
&z\in\Omega\setminus D_\epsilon,\;\;d(z,\partial{\Omega\setminus D_\epsilon})\geq C\epsilon^\frac{1}{2}\vert\ln{\epsilon}\vert^3,\end{split}\]
from \eqref{a?} with $a=a_-$ we conclude
\[\begin{split}
&z\in\Omega,\;\;d(z,\partial^+\Omega)\geq C^1\epsilon^\frac{1}{4}\vert\ln{\epsilon}\vert^\frac{3}{2},\\
&\Rightarrow\;\;\vert u(z)-a_-\vert<\delta.
\end{split}
\]
From this, proceeding as in the proof of \eqref{exp-refin} in the proof of Lemma \ref{nearone*} we derive
 the exponential estimate
\[\vert u(z)-a_-\vert\leq\delta e^{-\frac{k}{\epsilon}(d(z,\partial^+\Omega)-C^1\epsilon^\frac{1}{2}\vert\ln{\epsilon}\vert^3)},\;\;
z\in\Omega,\;d(z,\partial^+\Omega)\geq C^1\epsilon^\frac{1}{4}\vert\ln{\epsilon}\vert^\frac{3}{2}.\]
The estimate \eqref{exp-bound} follows from this, \eqref{bounds} and a suitable choice of $K$. The proof is complete.
\end{proof}
\begin{remark}
Theorem \ref{BL-exist} implies the existence of the pointwise limit
\[u_0=\lim_{\epsilon\rightarrow 0^+}u,\]
where
\[u_0=\left\{\begin{array}{l}
a_-,\;\;\text{on}\;\;\bar{\Omega}\setminus\partial^+\Omega,\\
a_+,\;\;\text{on}\;\;\partial^+\Omega.
\end{array}\right.\]
\end{remark}
\subsubsection{On the thickness of the boundary layer}
Fix a point $(\hat{x},\hat{y})\in\bar{\Omega}$ and consider the rescaled map
\begin{equation}
U^\epsilon(\xi,\eta)=u(\epsilon\xi+\hat{x},\epsilon\eta+\hat{y}),
\;\;(\epsilon\xi+\hat{x},\epsilon\eta+\hat{y})\in\bar{\Omega}.
\label{rescaled}
\end{equation}
The bound
\[\vert u\vert\leq M,\quad\;\Leftrightarrow\quad\;\vert U^\epsilon\vert\leq M,\]
the smoothness of $W$ and of $\partial\Omega$ and elliptic theory imply
\begin{equation}
\vert U^\epsilon\vert_{C^{2,\alpha}(\bar{\Omega}^\epsilon;\R^m)}\leq C.
\label{Ueps-bound}
\end{equation}
for some $\alpha\in(0,1)$ and some constant $C>0$.
Therefore, at least along a subsequence, there exists
\[U^0=\lim_{\epsilon\rightarrow0^+}U^\epsilon,\]
and the convergence is in $C^2$ sense in compact sets. Clearly the limit function $U^0$ depends
on the choice of the point $(\hat{x},\hat{y})\in\bar{\Omega}$. From Theorem \ref{BL-exist} we can
easily characterize $U^0$ for various choices of $(\hat{x},\hat{y})$. For instance if $(\hat{x},\hat{y})\in\bar{\Omega}\setminus\bar{\partial^+\Omega}$ we have
\[U^0=a_-,\;\;\text{for}\;\;(\xi,\eta)\in S,\]
where $S=\R^2$ if $(\hat{x},\hat{y})\in\Omega$ and $S$ is the half plane that contains $\Omega$ and is tangent to $\Omega$ at $(\hat{x},\hat{y})$ if $(\hat{x},\hat{y})\in\partial\Omega\setminus\bar{\partial^+\Omega}$. If $(\hat{x}^\epsilon,\hat{y})\in\partial\Omega$ has $\hat{x}^\epsilon\in(0,C\epsilon)$, $\hat{y}=0$ and $g(\hat{x},0)=b\neq a^\pm$ we expect that $U^0:S\rightarrow\R^m$ satisfies $\lim_{\xi\rightarrow\pm\infty}U^0(\xi,0)=a^\pm$ and that as $\eta\rightarrow+\infty$, $U^0(\cdot,\eta)$ converges exponentially to a translate of the heteroclinic connection
$\bar{u}$, see \cite{scha} and Chapter 9 in \cite{afs}. It remains to consider the case where $(\hat{x},\hat{y})\in\partial^+\Omega$. In this case, Theorem \ref{BL-exist} suggests that
\begin{equation}
U^0\equiv a_+,
\label{U00a}
\end{equation}
on the half plane $S$.
In spite of the fact that the estimate for the thickness of the boundary layer given by Theorem \ref{BL-exist} may  not be optimal as far as the power $\epsilon^\frac{1}{2}$, it is correct in indicating that the layer is thicker than $\mathrm{O}(\epsilon)$. In Theorem \ref{thick} below we establish \eqref{U00a} and also prove that the thickness is $\frac(\epsilon){\mathrm{o}(1)}$. This is compatible with the fact that there is no connecting orbit in the half plane.

The existence of the boundary layer is a higher dimensional effect and a new phenomenon. It is the result of a compromise between two competing minimization requirements: in the interior of $\partial^+\Omega$ to reduce energy the solution $U^\epsilon$ tries to behave like in the one dimensional case and push the layer in the interior while in a neighborhood of the extreme points of  $\partial^+\Omega$ minimization requires the solution to remain near $a_-$.
 \begin{theorem}
 \label{thick}
 \begin{enumerate}

Let $\hat{x}\in(0,l)$ be fixed.
Then
\item
For each $K>0$ it results
\[
\lim_{\epsilon\rightarrow 0^+}\sup_{y\leq K\epsilon}\vert u^{\epsilon}(\hat{x},y)-a_+\vert=0.\]
\item
\[\liminf_{\epsilon\rightarrow 0^+}\vert u^{\epsilon}(\hat{x},y^\epsilon)-a_+\vert>0\;\;\Rightarrow\;\;
\lim_{\epsilon\rightarrow 0^+}\frac{y^\epsilon}{\epsilon}=+\infty.\]
\item
\[\lim_{\epsilon\rightarrow 0^+}\epsilon\frac{\partial u^{\epsilon}}{\partial y}(\hat{x},0)=0.
\]
\end{enumerate}
\end{theorem}

 \begin{proof}
1. We begin with the scalar case $m=1$ that can be handled by simply combining known facts based on the Modica inequality which is not available in the vector case. Specifically we will make use of a result of Farina and Valdinoci [\cite{FaVa}, Theorem 1 (ii)] which in particular establishes the validity of the Modica inequality on the half space (as opposed to the whole $\R^n$).

 The map $U^\epsilon$ defined in \eqref{rescaled} for fixed $\hat{x}\in(0,l)$ and $\hat{y}=0$ satisfies
 \[\Delta U^\epsilon=W^\prime(U^\epsilon),\;\;(\epsilon\xi+\hat{x},\epsilon\eta)\in\Omega;
 \;\;\xi=\frac{x-\hat{x}}{\epsilon}, \eta=\frac{y}{\epsilon}.\]
 Since we are in the scalar case we can assume $a^\pm\in\R$ ,$a_-<a_+$. Then we have $a_-\leq U^\epsilon\leq a_+$ which follows from the boundary conditions on $\partial\Omega$ via the maximum principle. By linear elliptic theory
 \[U^\epsilon\in C^{2,\alpha}(\bar{\Omega}^\epsilon),\;\;
 \Omega^\epsilon=\{(\xi,\eta):(\epsilon\xi+\hat{x},\epsilon\eta)\in\Omega\},\]
 and so along a subsequence
 \[U^{\epsilon_j}\stackrel{C^2}{\rightarrow} U^0,\;\;\text{on compacts in }\,S=\R\times[0+\infty).\]
 The limit function $U^0$ satisfies
 \begin{equation}
 \left\{\begin{array}{l}
 \Delta U^0=W^\prime(U^0),\;\;\text{in}\,S=\R\times(0,+\infty),\\
 U^0(\xi,0)=a_+,\;\;\xi\in\R.
 \end{array}\right.
 \label{limitU}
 \end{equation}
 Setting $u=-U^0+a_+$, $-F^\prime(u)=W(a_+-u)$
 we can apply Theorem 1 (ii) in  \cite{FaVa} and conclude, via the positivity of $W$ and the specific bounds $ U^0\in[a_-,a_+]$, the Modica estimate
 \begin{equation}
 \frac{1}{2}\vert\nabla U^0\vert^2\leq W(U^0),\;\;\text{in}\;S.
 \label{Modic}
 \end{equation}
 As we mention below, classical estimates from linear elliptic theory imply
 \begin{equation}
 U^0\in C^{2,\alpha}(\bar{S}),
 \label{U-smooth}
 \end{equation}
 which we accept for the moment
 This extends the validity of \eqref{Modic} up to the boundary of $S$. Arguing now as in
 [\cite{Mo}, Theorem I] we set
 \[\phi(t)=U^0((\xi,0)+tn)+a_+,\;\;n\in\SF^1,\;t\in\R,\;(\xi,0)+tn\in S,\]
 and obtain via \eqref{Modic} that there is $\delta>0$ such that
 \[\vert\phi^\prime(t)\vert^2\leq C\vert\phi(t)\vert^2,\;\;\phi(0)=0,\;\vert t\vert\leq\delta,\]
 and conclude that $\phi(t)=0$ for $\vert t\vert<\delta$. This show that the set $\{(\xi,\eta)\in\bar{S}:U^0(\xi,\eta)=a_+\}$ which is nonempty and closed in $\bar{S}$ is also open and therefore coincides with $\bar{S}$, hence
 \begin{equation}
 U^0\equiv a_+,\;\;\text{on}\,\R\times[0,+\infty).
 \label{constant}
 \end{equation}
 From this Theorem \ref{thick} follows. To avoid repetition, we give details concerning  the various statements of the theorem, when dealing with the vector case.

 To prove \eqref{U-smooth} we note that the boundedness of $U^0$ and Theorem 8.29 in \cite{GT} imply a Holder estimate for $U^0$ up to the boundary of $S$. This and the smoothness of $W$ yields $W^\prime(U^0)\in C^\alpha(\bar{S})$, for some $\alpha\in(0,1)$. Next with the global Schauder estimate we have
 \[\|U^0\|_{C^{2,\alpha}(\bar{S})}\leq C_1\|W^\prime(U^0)\|_{C^\alpha(\bar{S})}+C_2\|U^0\|_{L^\infty(S)},\]
 and \eqref{U-smooth} follows.

 \begin{remark}
 The previous discussion of the scalar case can be slightly generalized. Indeed the argument developed to derive
 \eqref{constant} goes through in the same way also if we allow the point $\hat{x}\in(C\epsilon,l-C\epsilon)$ in the definition of $U^\epsilon$ to depend on $\epsilon$.
 \end{remark}

2. Now we move to the vector case $m>1$. The objective is to establish \eqref{constant} and the difficulty is due to the absence of the Modica estimate. We will utilize the upper bound \eqref{UB*} in Proposition \ref{UB2}, the refined lower bound \eqref{LB-ref} which is sharp as far as the power of $\epsilon$ is concerned as well as the zero order term, Gui's Hamiltonian identities (e.g. \cite{afs}, 3.4) and Theorem \ref{BL-exist} above.
 \begin{lemma}
 \label{Kin-bound}
\[\begin{split}
 &\int_0^l\int_0^h\vert u_x\vert^2dydx\leq C\vert\ln{\epsilon}\vert^3,\;\;\text{some constant}\;C>0,\\
 &\int_{-\frac{\hat{x}}{\epsilon}}^{\frac{l-\hat{x}}{\epsilon}}\int_0^{\frac{h}{\epsilon}}\vert U_\xi^\epsilon\vert^2 d\eta d\xi \leq C\vert\ln{\epsilon}\vert^3.
 \end{split}\]
 \end{lemma}
 \begin{proof}
 From \eqref{UB*} and \eqref{LB-ref} we have
 \[\frac{\epsilon}{2}\int_0^l\int_0^h\vert u_x\vert^2dydx\leq 2l\sigma+C_1\epsilon\vert\ln{\epsilon}\vert^3
 -\int_0^lJ^\epsilon_{(0,h)}(u(x,\cdot))dx\leq C_1\epsilon\vert\ln{\epsilon}\vert^3+\tilde{C}^\prime\epsilon,\]
 and the first inequality is established. The second inequality follows from the first by changing variables.
 \end{proof}
 From the boundary conditions we have $U^\epsilon(\xi,0)=a_+$ for $\xi\in(-\frac{\hat{x}}{\epsilon}+C,\frac{l-\hat{x}}{\epsilon}-C)$, so by linear elliptic theory (e.g.  Theorem 8.29 in \cite{GT})
 \begin{equation}
 \vert U^\epsilon(\xi,\eta)-U^\epsilon(\xi,0)\vert\leq\tilde{C}\vert\eta\vert^\alpha,
 \label{holder}
 \end{equation}
 for some constant $\tilde{C}>0$ and $\alpha\in(0,1)$.
 Thus we obtain in this range of $\eta$\,'s the estimate
  \begin{equation}
 \vert U^\epsilon(\xi,\eta)-a_+\vert\leq\tilde{C}\vert\eta\vert^\alpha,
 \label{holder1}
 \end{equation}
 \begin{lemma}
 \label{average0}
 There exist $\xi^\pm(\epsilon)$ such that
 \begin{equation}
 \begin{split}
 &\xi^-(\epsilon)\in(-\frac{\hat{x}}{\epsilon}+C, -\frac{\hat{x}}{\epsilon}+C+\epsilon^{-\frac{3}{4}}),\\
 &\xi^+(\epsilon)\in(\frac{l-\hat{x}}{\epsilon}-C-\epsilon^{-\frac{3}{4}},\frac{l-\hat{x}}{\epsilon}-C),\\
 &\int_0^{\frac{h}{\epsilon}}\vert U_\xi^\epsilon(\xi^\pm(\epsilon),\eta)\vert^2 d\eta\leq C\epsilon^\frac{3}{4}\vert\ln{\epsilon}\vert^3.
 \end{split}
 \label{average}
 \end{equation}
 \end{lemma}
 \begin{proof}
 Set $f(\xi)=\int_0^{\frac{h}{\epsilon}}\vert U_\xi^\epsilon(\xi,\eta)\vert^2 d\eta$, $\xi\in(-\frac{\hat{x}}{\epsilon}+C,\frac{l-\hat{x}}{\epsilon}-C)$. Then  Lemma \ref{Kin-bound} and the Mean Value Theorem imply, for each $p\in(-\frac{\hat{x}}{\epsilon}+C,\frac{l-\hat{x}}{\epsilon}-C-\epsilon^{-\frac{3}{4}})$
  \[\epsilon^{-\frac{3}{4}}f(\xi(\epsilon))=\int_{p}^{p+\epsilon^{-\frac{3}{4}}}f(\xi)d\xi\leq C\vert\ln{\epsilon}\vert^3,\;\;\text{some}\;\xi(\epsilon)\in(p,p+\epsilon^{-\frac{3}{4}}),\]
 \end{proof}

 We now utilize the Hamiltonian Identity in the rectangle $[\xi^-(\epsilon),\xi^+(\epsilon)]\times[0,\frac{h}{2\epsilon}]$\; arguing as in the proof of Theorem 3.2 in \cite{afs}. Taking into account the boundary condition $U^\epsilon(\xi,0)=a_+$, we have
 \begin{equation}
 \begin{split}
 &\int_{\xi^-(\epsilon)}^{\xi^+(\epsilon)}[\frac{1}{2}(\vert U_\xi^\epsilon\vert^2-\vert U_\eta^\epsilon\vert^2)
 +W(U^\epsilon)]\vert_{\eta=\frac{h}{2\epsilon}}d\xi+\frac{1}{2}\int_{\xi^-(\epsilon)}^{\xi^+(\epsilon)}\vert U_\eta^\epsilon(\xi,0)\vert^2d\xi\\
 &=\int_0^{\frac{h}{2\epsilon}}U_\xi^\epsilon\cdot U_\eta^\epsilon\vert_{\xi=\xi^-(\epsilon)}d\eta-
 \int_0^{\frac{h}{2\epsilon}}U_\xi^\epsilon\cdot U_\eta^\epsilon\vert_{\xi=\xi^+(\epsilon)}d\eta.
 \end{split}
 \label{Hamilton}
 \end{equation}

 \begin{lemma}
 \label{UxiUeta}
 \begin{equation}
 \begin{split}
 &\lim_{\epsilon\rightarrow 0^+}\int_{\xi^-(\epsilon)}^{\xi^+(\epsilon)}[\frac{1}{2}(\vert U_\xi^\epsilon\vert^2-\vert U_\eta^\epsilon\vert^2)
 +W(U^\epsilon)]\vert_{\eta=\frac{h}{2\epsilon}}d\xi=0,\\
 &\lim_{\epsilon\rightarrow 0^+}\int_0^{\frac{h}{2\epsilon}}\vert U_\xi^\epsilon\cdot U_\eta^\epsilon\vert\left\vert_{\xi=\xi^\pm(\epsilon)}\right. d\eta=0.
 \end{split}
 \label{UxiUeta1}
 \end{equation}
  \end{lemma}
 \begin{proof}
 From $\Delta U^\epsilon=W_u(U^\epsilon)=W_u(U^\epsilon)-W_u(a_-)$ and Theorem \ref{BL-exist}, via a local $L^p$-estimate, we obtain
 \begin{equation}
 \vert\nabla U^\epsilon\vert\leq Ce^{-k(\min\{\eta,\frac{h}{\epsilon}-\eta\}-C_1\frac{\vert\ln{\epsilon}\vert^3}{\epsilon^\frac{1}{2}})},\;\;\text{in}\;
 [-\frac{\hat{x}}{\epsilon},\frac{l-\hat{x}}{\epsilon}]\times[0,\frac{h}{\epsilon}],
 \label{Ueta}
 \end{equation}
 where $C>0$ here and below is a constant possibly different from line to line. This with $\eta=\frac{h}{2\epsilon}$ and the smallness of $W(U^\epsilon(\xi,\frac{h}{2\epsilon}))$, $\xi\in(-\frac{\hat{x}}{\epsilon},\frac{l-\hat{x}}{\epsilon})$ that follows by \eqref{exp-bound}, proves \eqref{UxiUeta1}$_1$.

 To complete the proof we write
\[\int_0^{\frac{h}{2\epsilon}}\vert U_\xi^\epsilon\cdot U_\eta^\epsilon\vert\left\vert_{\xi=\xi^\pm(\epsilon)}\right. d\eta=
 \int_0^{C_1\frac{\vert\ln{\epsilon}\vert^3}{\epsilon^\frac{1}{2}}}\vert U_\xi^\epsilon\cdot U_\eta^\epsilon\vert\left\vert_{\xi=\xi^\pm(\epsilon)}\right. d\eta+
 \int_{{C_1\frac{\vert\ln{\epsilon}\vert^3}{\epsilon^\frac{1}{2}}}}^{\frac{h}{2\epsilon}}\vert U_\xi^\epsilon\cdot U_\eta^\epsilon\vert\left\vert_{\xi=\xi^\pm(\epsilon)}\right. d\eta= I+II.\]
 We estimate each term separately. From \eqref{Ueps-bound} and \eqref{average}$_2$ and \eqref{Ueta} we have
 \begin{equation}
 I\leq C(\int_0^{C_1\frac{\vert\ln{\epsilon}\vert^3}{\epsilon^\frac{1}{2}}}\vert U_\xi^\epsilon(\xi^\pm(\epsilon),\eta)\vert^2 d\eta)^\frac{1}{2}(\frac{\vert\ln{\epsilon}\vert^3}{\epsilon^\frac{1}{2}})^\frac{1}{2}
 \leq C(\epsilon^\frac{3}{4}\vert\ln{\epsilon}\vert^3)^\frac{1}{2}
 (\frac{\vert\ln{\epsilon}\vert^3}{\epsilon^\frac{1}{2}})^\frac{1}{2}
 \leq C\epsilon^\frac{1}{8}\vert\ln{\epsilon}\vert^3,
 \label{I}
 \end{equation}
  \begin{equation}
  \begin{split}
  & II\leq(\int_0^{\frac{h}{\epsilon}}\vert U_\xi^\epsilon(\xi^\pm(\epsilon),\eta)\vert^2 d\eta)^\frac{1}{2}
  (\int_{C_1\frac{\vert\ln{\epsilon}\vert^3}{\epsilon^\frac{1}{2}}}^{\frac{h}{2\epsilon}}\vert U_\eta^\epsilon(\xi^\pm(\epsilon),\eta)\vert^2 d\eta)^\frac{1}{2}\\
  &\leq C(\epsilon^\frac{3}{4}\vert\ln{\epsilon}\vert^3)^\frac{1}{2}
  (\int_{C_1\frac{\vert\ln{\epsilon}\vert^3}{\epsilon^\frac{1}{2}}}^{\frac{h}{2\epsilon}}
  e^{-2k(\eta-C_1\frac{\vert\ln{\epsilon}\vert^3}{\epsilon^\frac{1}{2}})}d\eta)^\frac{1}{2}
  \leq C(\epsilon^\frac{3}{4}\vert\ln{\epsilon}\vert^3)^\frac{1}{2}.
  \end{split}
 \label{II}
 \end{equation}
 \end{proof}
 From Lemma \ref{UxiUeta} and \eqref{Hamilton} we obtain

 \begin{equation}
 \lim_{\epsilon\rightarrow 0^+}\int_{\xi^-(\epsilon)}^{\xi^+(\epsilon)}\vert U_\eta^\epsilon(\xi,0)\vert^2d\xi=0,
 \label{Ueta-small}
 \end{equation}
 and since, along a sequence $\{\epsilon_j\}$, $U^{\epsilon_j}$ converges locally in $\bar{S}=\R\times[0,+\infty)$, to $U^0$ in the $C^2$ sense, we conclude
 \begin{equation}
 \int_\R\vert U_\eta^0(\xi,0)\vert^2d\xi=0,
 \label{Ueta0-small}
 \end{equation}
 hence
  \begin{equation}
 U_\eta^0(\xi,0)=0,\;\;\xi\in\R.
 \label{Ueta0-loc}
 \end{equation}
 Passing to the limit as $\epsilon\rightarrow 0^+$ in \eqref{holder1} we also have
  \begin{equation}
 U^0(\xi,0)=0,\;\;\xi\in\R.
 \label{U0-loc}
 \end{equation}
 A classical argument based on \eqref{Ueta0-small} and \eqref{U0-loc} shows that the map $\tilde{U}:\R^2\rightarrow\R^m$
 \[\tilde{U}=\left\{\begin{array}{l}
 a_+,\;\;\text{on}\;\R\times(-\infty,0),\\
 U^0,\;\;\text{on}\;\R\times[0,+\infty),
 \end{array}\right.\]
 is a $W^{1,2}$ solution of $\Delta U=W_u(U)$. Obviously the same is true for the map identically equal to $a_+$. This and a unique continuation theorem in \cite{gl} imply $\tilde{U}\equiv a_+$ and therefore that
  \eqref{constant} holds also in the vector case.

  We are now in the position to conclude the proof. Assume that (i) does not hold. Then there exists $\bar{K}>0$, $\delta>0$ and sequences $\{\epsilon_j\}$, $\epsilon_j\rightarrow 0^+$ as $j\rightarrow+\infty$, $\{y_j\}$ such that
  \[\vert u^{\epsilon_j}(\hat{x},y_j)-a_+\vert\geq\delta,\;\;y_j\leq\bar{K}\epsilon_j,\;j=1,2,\ldots
 \]
  This is equivalent to
  \[\vert U^{\epsilon_j}(0,\eta_j)-a_+\vert\geq\delta,\;\;\eta_j\leq\bar{K},\;j=1,2,\ldots
 \]
 which, since $U^{\epsilon}$ converges to $U^0$ uniformly in compacts, contradicts \eqref{constant}. This contradiction proves (i). Statement (ii) is an straightforward  consequence of (i). Finally if (iii) does not hold we have
 \[\epsilon_ju_y^{\epsilon_j}(\hat{x},0)\geq\delta,\;\;j=1,2,\ldots\]
 for some sequence $\{\epsilon_j\}$, $\epsilon_j\rightarrow 0^+$ as $j\rightarrow+\infty$ and therefore
 \[U_\eta^{\epsilon_j}(0,0)\geq\delta,\;\;j=1,2,\ldots\]
 This again contradicts \eqref{constant} since, along a subsequence, $U_\eta^{\epsilon_j}$ converges in the $C^2$ sense, uniformly in compacts, to $U^0$. The proof of Theorem \ref{thick} is complete.

 \end{proof}

\subsection{$\frac{h}{l}<1$, The internal layer case}
In this section we analyze in detail the structure of the minimizers of problem \eqref{problem2} under the assumption
\begin{equation}
l>h.
\label{l>h}
\end{equation}
We will show that for small $\epsilon>0$, aside from two interior layers aligned with $\{0\}\times[0,h]$ and $\{l\}\times[0,h]$ across which minimizers switch from $a_-$ to $a_+$ and back from $a_+$ to $a_-$, minimizers stay close to $a_+$ inside the rectangle $(0,l)\times[0,h]$ and to $a_-$ outside. More precisely we have
\begin{theorem}
\label{TThh}
There are constants $k,K>0$ and $C>0$ independent of $\epsilon\in(0,\epsilon_0]$ for some $\epsilon_0>0$ such that
\begin{equation}
\left.\begin{array}{l}
\vert u(z)-a_-\vert\leq Ke^{-\frac{k}{\epsilon}(d(z,R)-C\epsilon^\frac{1}{4})^+},\;\;z\in \Omega\setminus R,\\\\
\vert u(z)-a_+\vert\leq Ke^{-\frac{k}{\epsilon}(d(z,\Omega\setminus R)-C\epsilon^\frac{1}{4})^+},\;\;
z\in R,
\end{array}\right.
\label{exp-TTh}
\end{equation}
where $R=(0,l)\times(0,h)$.
\end{theorem}
The proof of Theorem \ref{TThh} is quite elaborate and is based on the derivation of tight lower$/$upper energy bounds and on the study of the one-dimensional problem in Section \ref{ONE}.
\subsubsection{The upper bound.}
\begin{proposition}
\label{UUBB}
There exists $C>0$ independent of $\epsilon$ such that, if $u$ is a minimizer of \eqref{problem2}, then
\[J_\Omega^\epsilon(u)\leq 2\sigma h+C\epsilon.\]
\end{proposition}
\begin{proof}
1.
From \eqref{Omega} it follows that there is a number $\delta>0$ and smooth functions $p,q:[-\delta,l+\delta]\rightarrow\R$ such that
\[\begin{split}
&\Omega\cap[-\delta,l+\delta]\times[-\delta,h+\delta]=S_\delta,\\
& S_\delta=\{(x,y):y\in(p(x),q(x)),\,x\in[-\delta,l+\delta]\}
,
\end{split}\]

\begin{equation}
\begin{split}
& p(x)=0,\;\;x\in[0,l];\;q(x)=h,\;\;x\in[0,l]\\
& p^\prime(0)=p^\prime(l)=q^\prime(0)=q^\prime(l)=0,\\
& p^{\prime\prime}(0)=p^{\prime\prime}(l)=q^{\prime\prime}(0)=q^{\prime\prime}(l)=0.
\end{split}
\label{p e q}
\end{equation}

We define a test function $u_{\mathrm{test}}$ in several steps. We let $\delta_\epsilon>0$ be a small number to be chosen later and set
\begin{equation}
\begin{split}
&u_{\mathrm{test}}(z)=a_-,\;\;z\in\overline{\Omega\setminus S}_{\delta_\epsilon},\\
&u_{\mathrm{test}}(z)=a_+,\;\;z\in[\delta_\epsilon,l-\delta_\epsilon]\times[0,h].
\end{split}
\label{away-layer}
\end{equation}\vskip.2cm
2. To complete the definition of $u_{\mathrm{test}}$ we observe that $S_{\delta_\epsilon}$ is mapped onto the rectangle $R_{\delta_\epsilon}=[-\delta_\epsilon,l+\delta_\epsilon]\times[0,h]$ via the map $z=(x,y)\rightarrow\zeta=(\xi,\eta)$ defined by
\begin{equation}
\begin{split}
&\xi=x,\;\;x\in(-\delta_\epsilon,l+\delta_\epsilon),\\
&\eta=\frac{h(y-p(x))}{q(x)-p(x)},\;\;y\in(p(x),q(x)),\\
&\text{with inverse}\\
&x=\xi,\;\;\xi\in(-\delta_\epsilon,l+\delta_\epsilon),\\
&y=\frac{\eta(q(\xi)-p(\xi))}{h}+p(\xi),\;\;\eta\in(0,h).
\end{split}
\label{diffeo}
\end{equation}

 From \eqref{p e q}, if, as we do, $\delta_\epsilon$ is chosen sufficiently small, say $\delta_\epsilon\leq\epsilon^\frac{1}{2}$, we have that $p, q-h$, $p^\prime, q^\prime$ are $\mathrm{O}(\epsilon)$
  for $x\in[-\delta_\epsilon,l+\delta_\epsilon]$. It follows
\begin{equation}
\begin{split}
&\frac{\partial\zeta}{\partial z}(z)=I+\mathrm{O}(\epsilon),\\
&\mathrm{det}\frac{\partial z}{\partial\zeta}(\zeta)=1+\mathrm{O}(\epsilon),
\end{split}
\label{Jacob}
\end{equation}
where $\frac{\partial\zeta}{\partial z}$ is the jacobian matrix of $\zeta$, $\mathrm{det}M$ the determinant of the matrix $M$ and $z:S_{\delta_\epsilon}\rightarrow R_{\delta_\epsilon}$ the inverse of $\zeta$. From \eqref{Jacob} it follows that, if $v:R_{\delta_\epsilon}\rightarrow\R^m$ is a smooth bounded map and if we  set

\[u_{\mathrm{test}}(z)=v(\zeta(z)),\;\;z\in S_{\delta_\epsilon},\]
then it results
\begin{equation}
\begin{split}
&\int_{S_{\delta_\epsilon}}\vert\nabla u_{\mathrm{test}}(z)\vert^2 dz
\leq\int_{R_{\delta_\epsilon}}\vert\frac{\partial\zeta}{\partial z}(z(\zeta))\vert^2\vert\nabla v(\zeta)\vert^2
\vert\mathrm{det}\frac{\partial z}{\partial\zeta}(\zeta)\vert d\zeta,\\
&\leq(1+\mathrm{o}(\epsilon))\int_{R_{\delta_\epsilon}}\vert\nabla v(\zeta)\vert^2d\zeta,\\
&\int_{S_{\delta_\epsilon}}W(u_{\mathrm{test}}(z))dz=\int_{R_{\delta_\epsilon}}W(v(\zeta))
\vert\mathrm{det}\frac{\partial z}{\partial\zeta}(\zeta)\vert d\zeta,\\
&\leq(1+\mathrm{o}(\epsilon))\int_{R_{\delta_\epsilon}}W(v(\zeta))d\zeta.
\end{split}
\label{S-R}
\end{equation}
3. The estimates \eqref{S-R} show that we can work on $R_{\delta_\epsilon}$. Obviously we need to impose $v=a_+$ on  $[\delta_\epsilon,l-\delta_\epsilon]\times[0,h]$ we obviously take $v=a_+$ and to ensure that $u_{\mathrm{test}}$ satisfies the boundary conditions in the minimization problem \eqref{problem2} we must require
\[v(\zeta)=\tilde{g}_\epsilon(\zeta)=:g_\epsilon(z(\zeta)),
\quad\zeta\in[-\delta_\epsilon,l+\delta_\epsilon]\times\{0,h\},\]
where $g_\epsilon$ is as in \eqref{g-def}. More explicitly we have

\begin{equation}
\begin{split}
&\tilde{g}=a_-,\;\;(\xi,\eta)\in[-\delta_\epsilon,0]\cup[l,l+\delta_\epsilon]\times\{0,h\},\\
&\tilde{g}=g_\epsilon,\;\;(\xi,\eta)\in[0,C_0\epsilon]\cup[l-C_0\epsilon,l]\times\{0,h\},\\
&\tilde{g}=a_+,\;\;(\xi,\eta)\in[C_0\epsilon,\delta_\epsilon]\cup[l-\delta_\epsilon,l-C_0\epsilon]\times\{0,h\},
\end{split}
\label{g-BC}
\end{equation}
with $g_\epsilon$ and $C_0$ as in \eqref{g-def}.
\vskip.2cm
4. We complete the definition of $v$. We divide $[-\delta_\epsilon,\delta_\epsilon]\cup[l-\delta_\epsilon,l+\delta_\epsilon]\times[0,h]$ in parts denoted $A,B,C,D$ as sketched in Fig. \ref{klematis}. The definition of $v$ in regions labeled with the same letter is similar. Therefore we only define $v$ on one (marked by the superscript $*$) of the regions denoted with the same letter. The extension to the entire domain is then straightforward. We set $\delta_\epsilon=\epsilon+\epsilon\lambda_\epsilon$, with $\lambda_\epsilon=\vert\ln{\epsilon^n}\vert$ and $n\geq 1$ a sufficiently large number and define
\begin{equation}
\begin{split}
&v(\xi,\eta)=\bar{u}(\frac{\xi}{\epsilon}),\;\;(\xi,\eta)
\in A^*=[-\epsilon\lambda_\epsilon,\epsilon\lambda_\epsilon]\times[\epsilon,h-\epsilon],\\
&v(\xi,\eta)=\frac{\eta}{\epsilon}(\bar{u}(\frac{\xi}{\epsilon})-\tilde{g}(\xi,0))+\tilde{g}(\xi,0),\;\;(\xi,\eta)
\in B^*=[-\epsilon\lambda_\epsilon,\epsilon\lambda_\epsilon]\times[0,\epsilon],\\
&v(\xi,\eta)=(1-(\frac{\xi}{\epsilon}-\lambda_\epsilon))(\bar{u}(\lambda_\epsilon)-a_+)+
a_+,\\
&\;\;(\xi,\eta)
\in C^*=[\epsilon\lambda_\epsilon,\epsilon\lambda_\epsilon+\epsilon]\times[\epsilon,h-\epsilon],\\
&v(\xi,\eta)=(1-(\frac{\xi}{\epsilon}-\lambda_\epsilon))\frac{\eta}{\epsilon}(\bar{u}(\lambda_\epsilon)-a_+)+a_+
,\\
&\;\;(\xi,\eta)
\in D^*=[\epsilon\lambda_\epsilon,\epsilon\lambda_\epsilon+\epsilon]\times[0,\epsilon].
\end{split}
\label{v-def}
\end{equation}

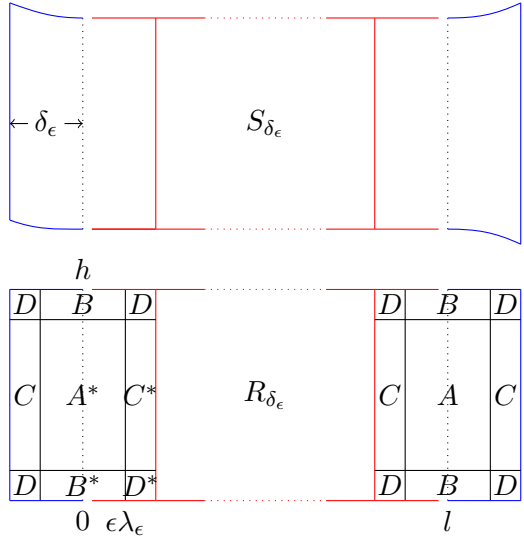
\begin{figure}
  \begin{center}
\begin{tikzpicture}[scale=.8]
\draw [] (-2.3,0) -- (-2.3,3.5);
\draw [red] (-1.8,0) -- (-1.8,3.5);
\node [] at (-2.05,.25) {$D^*$};

\draw [] (-3.7,0) -- (-3.7,3.5);
\draw [blue] (-4.2,0) -- (-4.2,3.5);
\node [] at (-3.95,.25) {$D$};

\draw [blue] (-4.2,0) -- (-3,0);
\draw [red] (-2.85,0) -- (-1.8,0);
\draw [] (-4.2,.5) -- (-1.8,.5);
\draw [] (-4.2,3) -- (-1.8,3);
\draw [blue] (-4.2,3.5) -- (-3,3.5);
\draw [red] (-2.85,3.5) -- (-1.8,3.5);

\draw [dotted] (-3,0) -- (-3,3.5);

\node [] at (-3.95,3.25) {$D$};
\node [] at (-3.,1.75) {$A^*$};
\node [] at (-3.,.25) {$B^*$};
\node [] at (-3.,3.25) {$B$};
\node [] at (-2.05,3.25) {$D$};
\node [] at (-2.05,1.75) {$C^*$};
\node [] at (-3.95,1.75) {$C$};

\node [below] at (-3,0) {$0$};
\node [above] at (-3,3.5) {$h$};

\node [below] at (-2.3,0) {$\epsilon\lambda_\epsilon$};


\draw [] (2.3,0) -- (2.3,3.5);
\draw [red] (1.8,0) -- (1.8,3.5);
\node [] at (2.05,.25) {$D$};

\draw [] (3.7,0) -- (3.7,3.5);
\draw [blue] (4.2,0) -- (4.2,3.5);
\node [] at (3.95,.25) {$D$};

\draw [red] (2.85,0) -- (1.8,0);
\draw [blue] (4.2,0) -- (3,0);
\draw [] (4.2,.5) -- (1.8,.5);
\draw [red] (2.85,3.5) -- (1.8,3.5);
\draw [blue] (4.2,3.5) -- (3,3.5);
\draw [] (4.2,3) -- (1.8,3);

\draw [dotted] (3,0) -- (3,3.5);

\node [] at (3.95,3.25) {$D$};
\node [] at (3.,1.75) {$A$};
\node [] at (3.,.25) {$B$};
\node [] at (3.,3.25) {$B$};
\node [] at (2.05,3.25) {$D$};
\node [] at (2.05,1.75) {$C$};
\node [] at (3.95,1.75) {$C$};

\draw [dotted,red] (-1,0) -- (1,0);
\draw [dotted,red] (-1,3.5) -- (1,3.5);
\draw [red] (-1.8,0) -- (-1,0);
\draw [red] (-1.8,3.5) -- (-1,3.5);
\draw [red] (1.8,0) -- (1,0);
\draw [red] (1.8,3.5) -- (1,3.5);

\node [] at (0,1.75) {$R_{\delta_\epsilon}$};

\node [below] at (3,0) {$l$};


\draw [dotted,red] (-1,4.5) -- (1,4.5);
\draw [dotted,red] (-1,8.) -- (1,8.);
\draw [red] (-1.8,4.5) -- (-1,4.5);
\draw [red] (-1.8,8.) -- (-1,8.);
\draw [red] (1.8,4.5) -- (1,4.5);
\draw [red] (1.8,8.) -- (1,8.);

\node [] at (0,6.25) {$S_{\delta_\epsilon}$};
\node [] at (-3.6,6.25) {$\delta_\epsilon$};
\draw [->] (-3.3,6.25)-- (-3,6.25);
\draw [->] (-3.9,6.25) -- (-4.2,6.25);

\draw [dotted] (3,4.5) -- (3,8);
\draw [dotted] (-3,4.5) -- (-3,8);

\draw [red] (-1.8,4.5) -- (-2.85,4.5);
\draw [red] (-1.8,4.5) -- (-2.85,4.5);
\draw [red] (-1.8,8) -- (-2.85,8);
\draw [red] (1.8,4.5) -- (2.85,4.5);
\draw [red] (1.8,8) -- (2.85,8);

\draw[blue] (3,4.5) to [out=0,in=155] (4.2,4.25);
\draw[blue] (3,8) to [out=0,in=205] (4.2,8.25);

\draw [blue] (4.2,4.25) -- (4.2,8.25);

\draw[blue] (-4.2,8.25) to [out=-20,in=180] (-3,8);
\draw[blue] (-4.2,4.65) to [out=-20,in=180] (-3,4.5);

\draw [blue] (-4.2,4.65) -- (-4.2,8.25);
\draw [red] (1.8,4.5) -- (1.8,8);
\draw [red] (-1.8,4.5) -- (-1.8,8);
\end{tikzpicture}
\end{center}
\caption{$R_{\delta_\epsilon}$, $S_{\delta_\epsilon}$ and the sets $A,B,C,D$.}
\label{klematis}
\end{figure}

 with these definitions we immediately have
 \begin{equation}
 J^\epsilon_{A^*}(v)\leq\sigma h.
 \label{Astar}
 \end{equation}
 To proceed with the estimates of $J_{B^*}(v)$ etc, we recall \eqref{kK} and that similar estimates are valid for the derivative of $\bar{u}$
. Moreover the fact that $v$ is bounded and $\mathbf{h}_1$ imply that there is a constant $C>0$ such that
 \begin{equation}
 W(v)\leq C\vert v-a^\pm\vert^2.
 \label{W-cv2}
 \end{equation}
 To estimate $J^\epsilon_{B^*}(v)$ we divide $B^*$ as $B^*=B^-\cup B^\prime\cup B^+$ where $B^-=[-\epsilon\lambda_\epsilon,0]\times[0,\epsilon]$, $B^\prime=[0,C_0\epsilon]\times[0,\epsilon]$ and
 $B^+=[C_0\epsilon,\epsilon\lambda_\epsilon]\times[0,\epsilon]$. From the properties of $g_\epsilon$ in \eqref{g-def} and the definition of $v$ it follows that $(\xi,\eta)\in B^\prime$ implies $\vert\nabla v\vert\leq\frac{C}{\epsilon}$ and $W(v)\leq C$. This and $\vert B^\prime\vert\leq C\epsilon^2$ imply
  \begin{equation}
 J^\epsilon_{B^\prime}(v)\leq C\epsilon.
 \label{Bprime}
 \end{equation}
 Next we estimate $J^\epsilon_{B^+}(v)$. Since, for $\xi\in[C_0\epsilon,\epsilon\lambda_\epsilon]$, it results $\tilde{g}(\xi,0)=a_+$ and therefore $v(\xi,\eta)=\frac{\eta}{\epsilon}(\bar{u}(\frac{\xi}{\epsilon})-a_+)+a_+$, from \eqref{kK} and \eqref{W-cv2} we obtain
 \[\begin{split}
 & J^\epsilon_{B^+}(v)\leq\frac{C}{\epsilon}\int_0^\epsilon\int_{C_0\epsilon}^{\epsilon\lambda_\epsilon}
 \Big(\vert\bar{u}(\frac{\xi}{\epsilon})-a_
 +\vert^2+\frac{\eta^2}{\epsilon^2}\vert\bar{u}^\prime(\frac{\xi}{\epsilon})\vert^2\Big)
 d\xi d\eta\\
 &\leq\frac{C}{\epsilon}\int_0^\epsilon\int_{C_0\epsilon}^{\epsilon\lambda_\epsilon}
 e^{-2k\frac{\xi}{\epsilon}}d\xi d\eta\leq C\epsilon.
 \end{split}\]
 In a similar way we also obtain $J^\epsilon_{B^-}(v)\leq C\epsilon$. From this and \eqref{Bprime} we conclude
 \begin{equation}
 J^\epsilon_{B^*}(v)\leq C\epsilon.
 \label{Bstar}
 \end{equation}
 On $C^*$ we have $\vert\frac{\partial v}{\partial\xi}\vert\leq\vert\bar{u}(\lambda_\epsilon)-a_+\vert$, $\frac{\partial v}{\partial\eta}=0$ and, from \eqref{W-cv2},  $W(v)\leq C\vert\bar{u}(\lambda_\epsilon)-a_+\vert^2$. This and \eqref{kK} imply
\begin{equation}
\begin{split}
&J^\epsilon_{C^*}(v)\leq\frac{C}{\epsilon}\vert\bar{u}(\lambda_\epsilon)-a_+\vert^2
\int_{\epsilon\lambda_\epsilon}^{\epsilon\lambda_\epsilon+\epsilon}d\xi
\\
&\leq Ce^{\ln{\epsilon^{2nk}}}\leq C\epsilon,
\end{split}
\label{Cstar}
\end{equation}
 provided
  $n$ is chosen sufficiently large.
 Finally, using as before \eqref{kK}, \eqref{W-cv2} and $\vert\bar{u}(\lambda_\epsilon)-a_+\vert^2\leq Ce^{\ln{\epsilon^{2nk}}}\leq C\epsilon$, it is seen that $J^\epsilon_{D^*}(v)$ if of higher order
 \begin{equation}
J^\epsilon_{D^*}(v)\leq C\epsilon^2.
\label{Dstar}
\end{equation}

From \eqref{Astar}, \eqref{Bstar}, \eqref{Cstar} and \eqref{Astar} and the similar estimates valid for $J_A(v)$ etc. we obtain
\[
J^\epsilon_{R_{\delta_\epsilon}}(v)\leq 2\sigma h+C\epsilon.
\]

This concludes the proof
\end{proof}
\subsubsection{The lower bound}
\begin{proposition}
 \label{LBex3}
 There exist $C_4>0$ and $\epsilon_0>0$ such that, if $u$ is a minimizer of \eqref{problem2}, then
 \[J^\epsilon_\Omega(u)\geq\int_\Omega\Big(\frac{\epsilon}{2}\vert\frac{\partial u}{\partial x}\vert^2+\frac{1}{\epsilon}W(u)\Big)dx dy\geq 2\sigma h-C_4\epsilon^{\frac{1}{2}},\;\;\epsilon\in(0,\epsilon_0].\]
\end{proposition}
\begin{remark}
The following refined lower bound is derived later, see \eqref{Ref-LB3} below
\[\int_\Omega\Big(\frac{\epsilon}{2}\vert\frac{\partial u}{\partial x}\vert^2+\frac{1}{\epsilon}W(u)\Big)dx dy
\geq 2\sigma h-\tilde{C}\epsilon.\]
As a Corollary to this and the upper bound above, one obtains the following key estimate
\[\int_\Omega\vert\frac{\partial u}{\partial y}\vert^2 dx dy\leq C.\]
\end{remark}
\begin{proof}
1. Let $\delta=\delta(\epsilon)>0$ a number to be chosen later that satisfies $\frac{\epsilon}{\delta^2}=\mathrm{o}(1)$, as $\epsilon\rightarrow 0$.
Set $\Omega_\delta=\{z\in\Omega:\min_{a\in A}\vert u(z)-a\vert<\delta\}$ and let $\Omega_\delta^c$ be the complement of $\Omega_\delta$ in $\Omega$.
 Then \eqref{outW} and Proposition \ref{UUBB} imply
\begin{equation}
\begin{split}
&\frac{c_W^2\delta^2}{2\epsilon}\vert\Omega_\delta^c\vert\leq\frac{1}{\epsilon}\int_\Omega W(u)dz\leq J_\Omega^\epsilon(u)\leq 2\sigma h+C\epsilon,\\
&\Rightarrow\;\;\vert\Omega_\delta^c\vert\leq C_5\frac{\epsilon}{\delta^2}.
\end{split}
\label{omegac-delta}
\end{equation}

2. For $x\in(0,l)$ set $\Sigma_x=\{x\}\times(0,h)$ and define
\[\begin{split}
& X=\{x\in(C_0\epsilon,l-C_0\epsilon):\mathcal{H}^1(\Sigma_x\cap\Omega_\delta)>h-\eta\frac{\epsilon}{\delta^2}\},\\
& X^c=\{x\in(C_0\epsilon,l-C_0\epsilon):\mathcal{H}^1(\Sigma_x\cap\Omega_\delta^c)\geq \eta\frac{\epsilon}{\delta^2}\},
\end{split}\]
 where $\eta>0$ is a constant to be selected later. From
\eqref{omegac-delta} we have
\[\eta\frac{\epsilon}{\delta^2}\mathcal{H}^1(X^c)\leq\int_{X^c}\mathcal{H}^1(\Sigma_x\cap\Omega_\delta^c)dx
\leq\vert\Omega_\delta^c\vert\leq C_5\frac{\epsilon}{\delta^2}.\]

It follows \[\mathcal{H}^1(X^c)\leq\frac{C_5}{\eta},\]
and therefore
\[\mathcal{H}^1(X)\geq l-2C_0\epsilon-\frac{C_5}{\eta}.\]

3. We divide the sections in $X$ in two parts $X_-$ and $X_+=X\setminus X_-$ where
\[\begin{split}
& X_-=\{x\in X: \exists\, a\in A\setminus\{a_+\}\,\text{and}\,z\in\Sigma_x\,\text{such that}\,\vert u(z)-a\vert<\delta\},\\
& X_+=\{x\in X: \Sigma_x\cap\Omega_\delta=\{z\in\Sigma_x:\vert u(z)-a_+\vert<\delta\}\}.
\end{split}\]
From the boundary condition at $x\in(0,l)$, $y=0,h$, Lemma \ref{lower-sigma} and Proposition \ref{UUBB} it follows
 that
\[\begin{split}
&(2\sigma-C_W\delta^2)\mathcal{H}^1(X_-)\leq 2\sigma h+C\epsilon,\\
&\Rightarrow\;\mathcal{H}^1(X_-)\leq\frac{h+\frac{C}{2\sigma}\epsilon}{1-\frac{C_W}{2\sigma}\delta^2}\leq h+\frac{C}{\sigma}\epsilon+\frac{C_Wh}{\sigma}\delta^2,
\end{split}\]
and in turn by Step 2. above
\[\mathcal{H}^1(X_+)\geq l-h-(2C_0+\frac{C}{\sigma})\epsilon-\frac{C_5}{\eta}-\frac{C_Wh}{\sigma}\delta^2.\]
Since $l>h$ this estimate shows that if $\epsilon>0$ and $\delta>0$ are sufficiently small and if $\eta>0$ is sufficiently large we have $\mathcal{H}^1(X_+)\geq\frac{1}{2}(l-h)$. Then from the characterization of $X_+$ and $X$ we have
\[x\in X_+\;\Rightarrow\;\mathcal{H}^1(\{z\in \Sigma_x:\vert u(z)-a_+\vert<\delta\})\geq h-\eta\frac{\epsilon}{\delta^2}.\]
This and Lemma \ref{lower-sigma} imply
\[\int_\Omega\Big(\frac{\epsilon}{2}\vert\frac{\partial u}{\partial x}\vert^2+\frac{1}{\epsilon}W(u)\Big)dx dy\geq(2\sigma-C_W\delta^2)(h-\eta\frac{\epsilon}{\delta^2})\geq 2\sigma h-(C_Wh+2\sigma\eta)\epsilon^\frac{1}{2},\]
where we have set $\delta=\epsilon^\frac{1}{4}$. The proof is complete.
\end{proof}

 For $z_1,z_2\in\R^2$ we denote $\mathrm{sg}(z_1,z_2)$ the open segment with end points $z_1$ and $z_2$. We use brackets for closed or half closed segments.
 For $y\in(0,h)$ let $\mathrm{sg}[z_y^1,z_y^2]$ the connected component of $\Sigma_y=\Omega\cap((-\infty,+\infty)\times\{y\})$ that contains $(0,l)\times\{y\}$.
 \begin{remark}
 \label{E-in-D1}
 From the proof of Proposition \ref{LBex3} there exist $0<x_1<x_2<l$ that satisfy $x_2-x_1\geq\frac{1}{2}(l-h)$ and are such that in
the lower bound in Proposition \ref{LBex3}, $\Omega$ can be replaced by the subset (see Figure \ref{D})
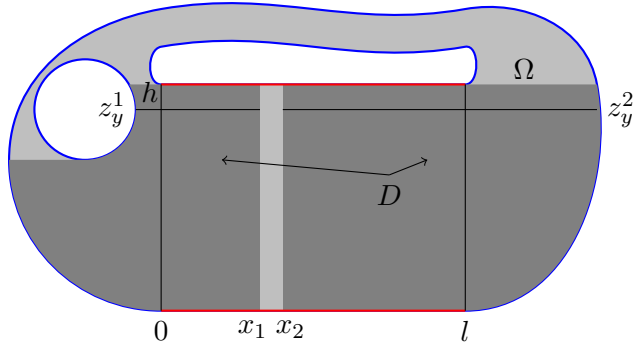
\begin{figure}
\begin{center}
\begin{tikzpicture}

\path [fill=lightgray] (0,0) -- (4,0) to [out=0,in=280] (5.732,3) to [out=100,in=10] (4,4)  to [out=190,in=10] (0,4)
to [out=190,in=90] (-2,2) to [out=270,in=180] (0,0);

\path [fill=white] (0,3) -- (4,3) to [out=0,in=10] (4,3.5) to [out=190,in=10] (0,3.5) to [out=190,in=180] (0,3);
\path [fill=white] (-1,2.6666) circle [radius=0.6666];;

\draw [blue, thick] (0,0) -- (4,0) to [out=0,in=280] (5.732,3) to [out=100,in=10] (4,4)  to [out=190,in=10] (0,4)
to [out=190,in=90] (-2,2) to [out=270,in=180] (0,0);

\draw [blue, thick] (0,3) -- (4,3) to [out=0,in=10] (4,3.5) to [out=190,in=10] (0,3.5) to [out=190,in=180] (0,3);
\draw [blue, thick]  (-1,2.6666) circle [radius=0.6666];

\path [fill=gray] (0,0) -- (4,0) to [out=0,in=280] (5.732,3) -- (-.42266,3) arc [radius=.6666, start angle=30, end angle=-90]--(-2,2) to [out=270,in=180] (0,0);

\path [fill=lightgray] (1.3,0) -- (1.6,0) -- (1.6,3) -- (1.3,3)-- (1.3,0);

\draw [red, thick] (0,0) -- (4,0);
\draw [red, thick] (0,3) -- (4,3);

\draw [thin](0,0)--(0,3);
\draw [thin](4,0)--(4,3);

\node[below] at (0,0) {$0$};
\node[below] at (4,0) {$l$};
\node[left] at (.1,2.9) {$h$};
\node[right] at (4.5,3.2) {$\Omega$};

\node[below] at (1.2,0) {$x_1$};
\node[below] at (1.7,0) {$x_2$};

\node[below] at (3,1.8) {$D$};
\draw [->](3,1.8)--(.8,2);
\draw [->](3,1.8)--(3.5,2);

\node[left] at (-.3334,2.6666) {$z_y^1$};
\node[right] at (5.732,2.6666) {$z_y^2$};

\draw (-.3334,2.6666)--(1.3,2.6666);
\draw (1.6,2.6666)--(5.732,2.6666);
\end{tikzpicture}

\end{center}
\caption{The set $D$ and the segments $\mathrm{sg}(z_y^1,(x_1,y))$ and $\mathrm{sg}((x_2,y),z_y^2)$.}
\label{D}
\end{figure}
\[D=\cup_{y\in(0,h)}\mathrm{sg}(z_y^1,(x_1,y))\cup\mathrm{sg}((x_2,y),z_y^2).\]
Here $D$ plays the same role as the set denoted again $D$ in the boundary layer case (see Remrk \ref{E-in-D})
\end{remark}
Arguing as in the proof of Theorem \ref{nearone}, we show that, in $\Omega\setminus D$, $u$ remains in a neighborhood of $A$. Actually, by adapting to the case at hand the argument in Steps 1. and 2. in the proof of Theorem \ref{nearone} we see that the existence of a point $z\in\Omega\setminus D$  that satisfies
\[d(z,\partial(\Omega\setminus D))\geq C\epsilon^\frac{1}{2},\;\;\text{ and}\;\; \min_{a\in A}\vert u(z)-a\vert\geq\delta,\]
for $C>0$ is sufficiently large, contradicts Proposition \ref{UUBB}. It follows
\begin{equation}
\begin{split}
& z\in\Omega\setminus D,\;\;d(z,\partial(\Omega\setminus D))\geq C\epsilon^\frac{1}{2}\\
&\Rightarrow\;\; \vert u(z)-a\vert\leq\delta,\;\;\text{for some}\;a\in A.
\end{split}
\label{first-dd}
\end{equation}
Set $Q=(x_1,x_2)\times(0,h)$ and $\Omega^*=\Omega\setminus(\overline{D\cup Q})$.
From \eqref{first-dd}, by means of the approach developed in Steps 2. and 3. in the proof of Lemma \ref{nearone*} we obtain that the continuity of $u$ and $u=g_\epsilon$ on $\partial\Omega$ imply

\begin{equation}
\begin{split}
& z\in Q,\;\;\min_i\vert x-x_i\vert\geq C\epsilon^\frac{1}{2}\\
&\Rightarrow\;\; \vert u(z)-a_+\vert\leq\delta,
\end{split}
\label{first-ddQ}
\end{equation}

and, if $\Omega^*$ is nonempty

\begin{equation}
\begin{split}
& z\in\Omega^*,\;\;d(z,D)\geq C\epsilon^\frac{1}{2}\\
&\Rightarrow\;\; \vert u(z)-a_-\vert\leq\delta.
\end{split}
\label{first-ddd}
\end{equation}

and then derive the exponential estimates
\begin{equation}
\begin{split}
&\vert u(z)-a_+\vert\leq K e^{-\frac{k}{\epsilon}(\min_i\vert x-x_i\vert-C\epsilon^\frac{1}{2})^+}, \;\;z\in Q,\\
&\vert u(z)-a_-\vert\leq K e^{-\frac{k}{\epsilon}(d(z,D)-C\epsilon^\frac{1}{2})^+}, \;\;z\in\Omega^*.
\end{split}
\label{need-est}
\end{equation}
The estimate \eqref{need-est}$_1$ implies that we can refine the lower bound in Proposition \ref{LBex3} to
\begin{equation}
\int_\Omega\Big(\frac{\epsilon}{2}\vert\frac{\partial u}{\partial x}\vert^2+\frac{1}{\epsilon}W(u)\Big)dx dy
\geq 2\sigma h-\tilde{C}\epsilon.
\label{Ref-LB3}
\end{equation}
Since $x_2-x_1\geq\frac{1}{2}(l-h)$, for $x=x_m=\frac{1}{2}(x_1+x_2)$ we have $\min_i\vert x_m-x_i\vert\geq\frac{1}{4}(l-h)$. This and \eqref{need-est}$_2$ yield, for $\epsilon>0$ sufficiently small,
\begin{equation}\vert u(x_m,y)-a_+\vert\leq Ke^{-\frac{k}{\epsilon}(\frac{1}{4}(l-h)-C\epsilon^\frac{1}{2})}\leq\epsilon^\frac{1}{2}.
\label{Along}
\end{equation}

Therefore Lemma \ref{lower-sigma} implies
\begin{equation}
J^\epsilon_{(x_y^1,x_m)}(u(\cdot,y)),\;\; J^\epsilon_{(x_m,x_y^2)}(u(\cdot,y))>\sigma-\frac{1}{2}C_W\epsilon,\;\;y\in(0,h),
\label{LB--l}
\end{equation}
where $x_y^i$ is defined by $(x_y^i,y)=z_y^i$, $i=1,2$,
  and \eqref{Ref-LB3} follows.

\subsubsection{The structure inside $D$ and the proof of Theorem \ref{TThh}.}
In this Section we parallel the reasoning developed in Section \ref{Structure} to analyze the fine structure of a minimizer $u$ inside the set $D$ and on the basis of this information we prove Theorem \ref{TThh}.
The estimates in \eqref{LB--l} suggest that $u(\cdot,y)\vert_{(x_y^1,x_m)}$ should be a perturbation of a translation of the heteroclinic $\bar{u}$ and similarly that $u(\cdot,y)\vert_{(x_m,x_y^2)}$ should be a perturbation of the map $s\rightarrow\bar{u}(-s)$ that connects $a_+$ to $a_-$.
By applying Lemma \ref{min-0,lambda} and Lemma \ref{winVc} to the restrictions of $u$ to the segments $\mathrm{sg}(z_y^1,(x_m,y))$ and $\mathrm{sg}((x_m,y),z_y^2)$ we show that this is true for most $y\in(0,h)$.
When dealing with $\mathrm{sg}(z_y^1,(x_m,y))$ the interval $(0,\lambda)$ in \eqref{problem*2} should be replaced by the interval $(x_y^1,x_m)$ where $x_y^1$ is defined by $(x_y^1,y)=z_y^1$, and the intervals $[0,s^{-,w}]$ and $
[s^{+,w},\lambda]$ with the intervals $[x_y^1,x_y^1+s^{-,w}]$ and $[x_y^1+s^{+,w},x_m]$.

Let $Y=Y_1\cup Y_2\subset(0,h)$ be defined by
\[\begin{split}
& Y_1=\{y\in(0,h):u(\cdot,y)\vert_{(x_y^1,x_m)}\in\mathscr{V}^c\},\\
& Y_2=\{y\in(0,h):u(\cdot,y)\vert_{(x_m,x_y^2)}\in\mathscr{V}^c\}.
\end{split}\]

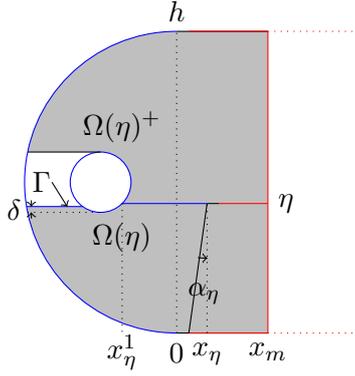
\begin{figure}
  \begin{center}
\begin{tikzpicture}[scale=.8]

\path[fill=lightgray] (0,-2.5)--(1.5,-2.5)--(1.5,-.3535)--(-.8964,-.3535)
arc [radius=0.5, start angle=-45, end angle=-126.136]--(-2.4672,-.4035)
arc [radius=2.5, start angle=189.288, end angle=270]--(0,-2.5);

\path[fill=lightgray] (1.5,-.3535)--(-.8964,-.3535)
arc [radius=0.5, start angle=-45, end angle=90]--(-2.449,.5)
arc [radius=2.5, start angle=168.463, end angle=90]--(1.5,2.5)--(1.5,-.3535);

\draw [blue] (0,2.5) arc [radius=2.5, start angle=90, end angle= 270];
\draw [blue] (-1.25,0) circle [radius=.5];

\draw[](-2.449,.5)--(-1.25,.5);
\draw[dotted](-2.449,-.5)--(-1.25,-.5);
\draw[blue](-2.4672,-.4035)--(-1.5453,-.4035);

\draw[](0,2.5)--(1.5,2.5);
\draw[red](.2,2.5)--(1.5,2.5);
\draw[dotted,red](1.5,2.5)--(3,2.5);

\draw[](0,-2.5)--(1.5,-2.5);
\draw[red](.2,-2.5)--(1.5,-2.5);
\draw[dotted,red](1.5,-2.5)--(3,-2.5);

\draw[dotted](0,-2.5)--(0,2.5);
\draw[red](1.5,-2.5)--(1.5,2.5);

\draw [blue] (-.8964,-.3535)--(.5,-.3535);
\draw [] (.5,-.3535)--(.7,-.3535);
\draw [red] (.7,-.3535)--(1.5,-.3535);

\draw[](.2,-2.5)--(.5,-.3535);
\draw[dotted](.5,-2.5)--(.5,-.3535);

\draw[->] (-2.389,-.3)--(-2.389,-.4035);
\draw[->] (-2.389,-.6035)--(-2.389,-.5);
\node[left] at (-2.389,-.45)  {$\delta$};
\draw[->] (-2.05,0)--(-1.80,-.4035);
\node[left] at (-1.9,0)  {$\Gamma$};

\node[above] at (-.9,.5)  {$\Omega(\eta)^+$};
\node[below] at (-.9,-.5)  {$\Omega(\eta)$};
\draw[dotted](-.8964,-.3535)--(-.8964,-2.5);

\node[below] at (-.8964,-2.35)  {$x_\eta^1$};
\node[below] at (0,-2.5)  {$0$};
\node[below] at (.5,-2.5)  {$x_\eta$};
\node[below] at (1.5,-2.5)  {$x_m$};
\node[right] at (1.5,-.3535)  {$\eta$};
\node[above] at (0,2.5)  {$h$};

\draw[->] (.35,-1.25) arc [radius=1.298, start angle=263.108, end angle=270];
\node[below] at (.45,-1.5)  {$\alpha_\eta$};
\end{tikzpicture}
\end{center}
\caption{$\Omega(\eta)$, $\Omega(\eta)^+$ and $\Gamma$.}
\label{Reta}
\end{figure}

\begin{lemma}
\label{bad-setin}
There exists $C>0$ such that
\begin{equation}
\mathcal{H}^1(Y)\leq C\epsilon^\frac{1}{2}.
\label{Y-bound}
\end{equation}
\end{lemma}
\begin{proof}
From Lemma \ref{winVc}, \eqref{LB--l} and Proposition \ref{UUBB} we obtain
\[\sum_{i\in\{1,2\}}\Big((\sigma-\frac{1}{2}C_W\epsilon)(h-\mathcal{H}^1(Y_i))+
\mathcal{H}^1(Y_i)(\sigma+C_W\epsilon^\frac{1}{2})\Big)\leq2\sigma h+C\epsilon.\]
It follows $\mathcal{H}^1(Y)\leq\mathcal{H}^1(Y_1)+\mathcal{H}^1(Y_2)\leq C\epsilon^\frac{1}{2}$.
\end{proof}
For $\eta\in(0,\frac{h}{2}]$ we define the set $\Omega(\eta)$, the analogous of the set  $\Omega(\xi)$ in the proof of Theorem \ref{BL-exist}. A set $\Omega$ that satisfies \eqref{Omega} may have a rather complex structure and the same is true for the set $D$, see Figure \ref{D}. Therefore some care is needed for the definition of $\Omega(\eta)$. We start by setting
\[\begin{split}
&\tilde{\Omega}(\eta)=\cup_{y\in(0,\eta)}\mathrm{sg}(z_y^1,(x_m,y)),\\
&\Omega(\eta)^+=\cup_{y\in[\eta,h)}\mathrm{sg}(z_y^1,(x_m,y)).
\end{split}\]
Set $\Pi(\eta)=\Omega^c\cup\Omega(\eta)^+\cup\{z:x\geq x_m\}$ and define
\[\Omega(\eta)=\tilde{\Omega}_{\delta}(\eta)\setminus\Pi(\eta),\]
where $\Omega^c$ is the complement of $\Omega$ in $\R^2$, $\tilde{\Omega}_{\delta}(\eta)=\cup_{z\in\tilde{\Omega}(\eta)}B_{\delta}(z)$ and $\delta>0$ is a small number.
The plan is now to obtain a precise estimate of $J^\epsilon_{\Omega(\eta)}(u)$. To do this we collect information of the values of the minimizer $u$ on $\partial\Omega(\eta)$.
The boundary of $\Omega(\eta)$ includes $\mathrm{sg}(z_\eta^1,(x_m,\eta))$ and $\partial\Omega(\eta)\cap\partial\Omega$. We observe that if
\[\Gamma=\partial\Omega(\eta)\setminus\Big(\mathrm{sg}(z_\eta^1,(x_m,\eta))
\cup(\partial\Omega(\eta)\cap\partial\Omega)\Big)\neq\emptyset,\]
then, for each $z\in\Gamma$ we have $d(z,D)\geq\delta$ (cfr. Figure \ref{Reta}).

This, \eqref{need-est}$_1$ and $u=a_-$ on $\partial\Omega^-$, imply that, for $\epsilon>0$ sufficiently small it results
\[\vert u(z)-a_-\vert<\epsilon,\;\;z\in\Gamma\cup(\partial\Omega(\eta)\cap\partial\Omega^- ).\]
Assume now that $\eta\not\in Y$ then we have $u(\cdot,\eta)\vert_{(x_\eta^1,x_m)}\in\mathscr{W}^*$. This, \eqref{v-fat} and \eqref{sw-sw} imply that, if we set $x_\eta=s^{-,u(\cdot,\eta)\vert_{(x_\eta^1,x_m)}}-x_\eta^1$, we have
\[\begin{split}
&\vert u(x,\eta)-a_-\vert\leq K\epsilon^\frac{1}{4},\;\;x\in(x_\eta^1,x_\eta),\\
&\vert u(x,\eta)-a_+\vert\leq K\epsilon^\frac{1}{4},\;\;x\in(x_\eta+2C^*\epsilon^\frac{1}{2},x_m).
\end{split}\]
To complete the description of the boundary values of $u$ on $\partial\Omega(\eta)$ we recall \eqref{Along} which implies
\[\vert u(x_m,y)-a_+\vert\leq\epsilon^\frac{1}{2},\;\;y\in(0,\eta),\]
and $u=g_\epsilon$ on $\partial\Omega$ and in particular
\[u(x,0)=a_+, \;\;x\in(C_0\epsilon,x_m).\]
We proceed to estimate $J^\epsilon_{\Omega(\eta)}(u)$. We give details for the case $\eta\in(0,\frac{h}{2})$ and $x_\eta>0$, the other cases can be discussed in the same way with obvious modifications. We regard $\Omega(\eta)$ as the union of fibers orthogonal to the segment $\mathrm{sg}((C_0\epsilon,0),(x_\eta,\eta))$. We let $\Omega(\eta)^\prime$ the union of the fibers that have one of the extreme on $\mathrm{sg}((C_0\epsilon,0),(x_m,0))$ and $\Omega(\eta)^{\prime\prime}$ the union of the fibers that have one of the extreme on $\mathrm{sg}((x_m,0),(x_m,\eta))$. From the above discussion on the boundary values of $u$ on $\partial\Omega(\eta)$ and Lemma \ref{lower-sigma} we obtain
\begin{equation}
\begin{split}
&J^\epsilon_{\Omega(\eta)^\prime}(u)\geq(\sigma-C_WK^2\epsilon^\frac{1}{2})(x_m-C_0\epsilon)\sin{\alpha_\eta},\\
&J^\epsilon_{\Omega(\eta)^{\prime\prime}}(u)\geq(\sigma-C_WK^2\epsilon^\frac{1}{2})
\Big(\frac{\eta}{\cos{\alpha_\eta}}-(x_m-C_0\epsilon)\sin{\alpha_\eta}\Big),\\
&\Rightarrow\\
&J^\epsilon_{\Omega(\eta)}(u)\geq(\sigma-C_WK^2\epsilon\epsilon^\frac{1}{2})
\frac{\eta}{\cos{\alpha_\eta}},
\end{split}
\label{Omprime}
\end{equation}
where
\[\sin{\alpha_\eta}=\frac{x_\eta-C_0\epsilon}{\sqrt{(x_\eta-C_0\epsilon)^2+\eta^2}},\quad
\cos{\alpha_\eta}=\frac{\eta}{\sqrt{(x_\eta-C_0\epsilon)^2+\eta^2}}.\]
On the other hand, since by construction $\Omega(\eta)\cap\Omega(\eta)^+=\emptyset$, from \eqref{LB--l} we obtain
\[\begin{split}
&J_\Omega^\epsilon(u)\geq J^\epsilon_{\Omega(\eta)}(u)+J^\epsilon_{D\setminus\Omega(\eta)}(u),\\
&J^\epsilon_{D\setminus\Omega(\eta)}(u)\geq(2h-\eta)(\sigma-C\epsilon).
\end{split}\]
This \eqref{Omprime} and Lemma \ref{UUBB} yield
\[(\sigma-C_WK^2\epsilon^\frac{1}{2})
\frac{\eta}{\cos{\alpha_\eta}}+(2h-\eta)(\sigma-C\epsilon)\leq 2h\sigma+C\epsilon,\]
which, recalling also the expression of $\cos{\alpha_\eta}$, implies

\begin{equation}
\begin{split}
&\sqrt{(x_\eta-C_0\epsilon)^2+\eta^2}-\eta\leq C\epsilon^\frac{1}{2},\\
&\Rightarrow\quad x_\eta\leq C\epsilon^\frac{1}{4}.
\end{split}
\label{xeta-small}
\end{equation}
Set
\[D_\epsilon=\Big((-C\epsilon^\frac{1}{4},C\epsilon^\frac{1}{4})\cup(l-C\epsilon^\frac{1}{4},l+C\epsilon^\frac{1}{4})\Big)
\times\Big((0,h)\setminus Y\Big).\]
Then \eqref{JW*}, Lemma \ref{bad-setin} and \eqref{xeta-small} imply
\[\begin{split}
&J^\epsilon_{D_\epsilon}(u)\geq2(\sigma-C_W\epsilon^\frac{1}{2})(h-\mathcal{H}^1(Y)),\\
&\Rightarrow\quad J^\epsilon_{D_\epsilon}(u)\geq2\sigma h-C\epsilon^\frac{1}{2}.
\end{split}\]
 This shows that most of the energy of $u$ is concentrated in $D_\epsilon$ and by consequence allows to apply the arguments developed for \eqref{need-est} to establish the exponential estimates in Theorem \ref{TThh}. The proof of Theorem  \ref{TThh} is complete.

\end{document}